\documentclass[paper=a4, english, 10pt]{amsart}

\usepackage{latexsym}
\usepackage{cite}
\usepackage{amsfonts}
\usepackage{amsmath}
\usepackage{amsthm}
\usepackage{amssymb}
\usepackage{mathrsfs}
\usepackage{mathabx}
\usepackage{bbm}           
\usepackage{dsfont}
\usepackage{psfrag}
\usepackage{epstopdf}

\allowdisplaybreaks

\usepackage{wrapfig}

\usepackage{framed}

\usepackage{tikz}
\usepackage{color}
\usepackage{graphicx}

\usepackage{caption}
\usepackage{subcaption}

\usepackage{fancyvrb}
\usepackage{cases}
\usepackage{multicol}
\usepackage{xspace}

\usepackage[utf8]{inputenc}
\usepackage[T1]{fontenc}

\usepackage{listingsutf8}

\usepackage{textcomp}
\usepackage{lmodern}
\usepackage[a4paper,margin=2cm]{geometry}
\usepackage[english]{babel}

\usepackage{hyperref}

\newtheorem{theorem}{Theorem}[section]

\newtheorem{definition}[theorem]{Definition}

\newcommand{\N}{\mathbbm{N}}

\newcommand{\R}{\mathbbm{R}}

\newcommand{\Z}{\mathbbm{Z}}

\newcommand{\Lp}[1]{\mathbf{L}^{#1}}

\newcommand{\Lip}{\mathbf{Lip}}

\newcommand{\sgn}[1]{\mathrm{sign}\left(#1\right)}

\newcommand{\pt}{\partial_t}

\newcommand{\px}{\partial_x }

\renewcommand{\d}{{\,\rm{d}}}

\begin{document}

\title[Validation ACD]{On the implementation of a finite volumes scheme with monotone transmission conditions for scalar conservation laws on a star-shaped network}

%
%

\author[]{Sabrina Francesca Pellegrino}
\address{Dipartimento di Matematica\\
Universit\`a degli Studi di Bari\\
Via Orabona 4\\
70125 Bari\\
Italy}
\email{sabrina.pellegrino@uniba.it}
\thanks{The author acknowledges the support of the R\'egion Bourgogne Franche-Comt\'e, projet 2017-2020 ``Analyse math\'ematique et simulation num\'erique d'EDP issus de probl\`emes de contr\^ole et du trafic routier"  and of the Universit\'e de Bourgogne Franche-Comt\'e, projet Chrysalide 2017 ``Contr\^ole, analyse num\'erique et applications d'\'equations hyperboliques sur un r\'eseau". The author is member of the Gruppo Nazionale per l'Analisi Matematica, la Probabilit\`a e le loro Applicazioni (GNAMPA) of the Istituto Nazionale di Alta Matematica (INdAM)}
%


\begin{abstract}

In this paper we validate the implementation of the numerical scheme proposed in~\cite{ACD}. The validation is made by comparison with an explicit solution here obtained, and the solutions of Riemann problems for several networks. We then perform some simulations in order to qualitatively validate the model under consideration.

Such results represent also a first step for the validation of the finite volumes scheme introduced in~\cite{DSDPR}.
\end{abstract}

\maketitle

\section{Introduction}
We investigate from the numerical point of view the model developed in \cite{ACD} by Andreianov, Coclite and Donadello, called here ACD, which describes the evolution of traffic at a junction consisting of $m$ incoming and $n$ outgoing arcs. Incoming arcs are parametrized by $x\in\R_{-}$ and numbered by the index $i\in\mathsf{I} = \{1,\ldots,m\}$, while outgoing arcs are parametrized by $x\in\R_+$ and numbered by the index $j\in\mathsf{J} = \{m+1,\ldots, m+n\}$ in such a way that the junction is always located at  $x=0$. We denote the generic arc by $\Omega_h$, $h\in\mathsf{H} = \{1,\ldots,m+n\}$, and the network by $\Gamma=\Pi_{h\in\mathsf{H}}\Omega_h$.

We describe the evolution of traffic on each arc by the Lighthill-Whitham-Richards (LWR) model~\cite{LWR1,LWR2}, namely by a scalar conservation law of the form
\begin{equation}
\label{eq:basic}
  \rho_{h,t} + f_h(\rho_h)_x  = 0, \qquad\text{for } t >0,\,x\in\Omega_h,\,h\in\mathsf{H},
\end{equation}
where $\rho_h$ is the density and $f_h$ is the flux on the $h$-th arc. 
We assume that the arcs have a common maximal density $\rho_{\max}>0$ and the fluxes are all bell-shaped (unimodal), Lipschitz and non-degenerate nonlinear i.e. 
\begin{enumerate}
  \item[\textbf{(F)}] for all $h\in\mathsf{H}$, $f_h \in \Lip\left( [0,\rho_{\max}]; \R_+ \right)$ with $\|f_h'\|_{\infty}\leq L_h$, $f_h(0) = 0 = f_h(\rho_{\max})$, and there exists $\rho_{h,c} \in \left(0,\rho_{\max}\right)$ such that $f_h'(\rho)~(\rho_{h,c}-\rho)>0$ for a.e.~$\rho \in [0,\rho_{\max}]$,
	
  \item[\textbf{(NLD)}] for all $h\in\mathsf{H}$,  $f'_h$ is not constant on any non-trivial subinterval of $[0,\rho_{\max}]$.
\end{enumerate}
We augment \eqref{eq:basic} with the initial conditions
\begin{equation}
\label{eq:init}
\rho_h(0,x)=\rho_{h,0}(x), \quad x\in\Omega_h,
\end{equation}
where $\rho_{h,0}\in\Lp\infty(\Omega_h ; [0,\rho_{\max}])$, $h\in\mathsf{H}$.
We also impose the conservation of the total density at the junction, i.e. for a.e. $t\in\R_+$
\begin{equation}
\label{conservation}
  \sum_{i\in\mathsf{I}} f_i\left(\rho_i(t,0^-)\right)  = \sum_{j\in\mathsf{J}} f_j\left(\rho_j(t,0^+)\right).
\end{equation}
Notice that the previous equation makes sense as the assumption {\bf (F)} ensures the existence of strong traces~\cite{Panov, Vasseur}.

In~\cite{ACD} the authors prove the well-posedness of solutions obtained as vanishing viscosity limits for the Cauchy problem~\eqref{eq:basic}-\eqref{eq:init}. Their result relies upon the explicit characterization of the class of admissible weak solution at the junction in terms of vanishing viscosity germ, see~\cite{AC,AKRvv,AKRL1}. It represents a generalized study of the model investigated in~\cite{CocliteGaravello2010}, in which the authors establish the existence of weak solutions as limit of vanishing viscosity approximations. Such results are relevant in the perspective of a theoretical analysis of PDEs on networks, in particular, they allow the extension of the analogy between vanishing viscosity and numerical scheme both in the network case. Furthermore, their analysis is applicable to general junction solvers enjoying enjoying the order-preservation property, see for instance~\cite{DSRDB}.

The aim of this paper is to validate the implementation of the numerical scheme proposed in \cite{ACD} by comparison with an explicit solution here computed and with the solutions of Riemann problems both for merge, divide and 2-2 cases. Moreover, we show the consistency of the scheme with respect to the case of a network with no discontinuity at the junction, namely taking the same flux on each arc, and finally a convergence analysis is also performed. These results are used in the validation of the finite volumes scheme with point constraints at the junction introduced in~\cite{DSDPR}.


The paper is organized as follows. In Section~\ref{sec:wellp} we briefly recall the main theoretical results for ACD. In Section~\ref{sec:explicit} we compute an explicit solution for the problem in the case of a merge consisting of two incoming and one outgoing arcs. In Section~\ref{sec:volume-scheme} we present the numerical scheme. Finally, Section~\ref{sec:validation} is devoted to the validation of the scheme. 


\section{Well-posedness of ACD in the frame of admissible solutions}
\label{sec:wellp}
The well-posedness for the general Cauchy problem~\eqref{eq:basic}-~\eqref{eq:init}-~\eqref{conservation} is established in \cite{ACD} in the frame of admissible solution.

We recall some definitions.
\begin{definition}
\label{def:weak-sol}
A function $\rho_h\in\Lp\infty(\R_+\times\Omega_h;[0,\rho_{\max}])$, $h\in\mathsf{H}$, is a {\em weak solution} of~\eqref{eq:IBVP} if
\begin{itemize}
\item 
for every $k\in[0,\rho_{\max}]$ and every nonnegative test function $\phi\in\mathcal{C}^{\infty}(\R\times\Omega_h;\R)$ with compact support
\begin{equation*}
\int_0^\infty\int_{\Omega_h}\Bigl(|\rho_h-k|\pt\phi + \sgn{\rho_h-k}\left(f_h(\rho_h)-f_h(k)\right)\px\phi\Bigr)\,\d x\d t + \int_{\Omega_h}\left|\rho_{h,0}(x)-k\right|\phi(0,x)\,\d x \ge0;
\end{equation*}
\item for a.e. $t>0$, it holds
\begin{equation}
\label{eq:cons}
\sum_{i\in\mathsf{I}} f_i(\rho_i(t,0-)) = \sum_{j\in\mathsf{J}} f_j(\rho_j(t,0+)).
\end{equation}
\end{itemize}
\end{definition}

We remind the formulation of the Bardos-LeRoux-N\'ed\'elec boundary condition for conservation laws in terms of the Godunov numerical flux (see~\cite{Bardos, LeFloch}), which will be useful for the definition of admissible solution at the junction.

\begin{figure}[ht]
\centering{
\begin{psfrags}
      \psfrag{G}[t,c]{$G(a,\cdot)$}
      \psfrag{a}[c,B]{$a$}
      \psfrag{b}[c,B]{$b$}
      \psfrag{r}[c,B]{$\rho_{\max}$}
      \psfrag{c}[c,B]{$\rho_c$}
\includegraphics[width=.23\textwidth]{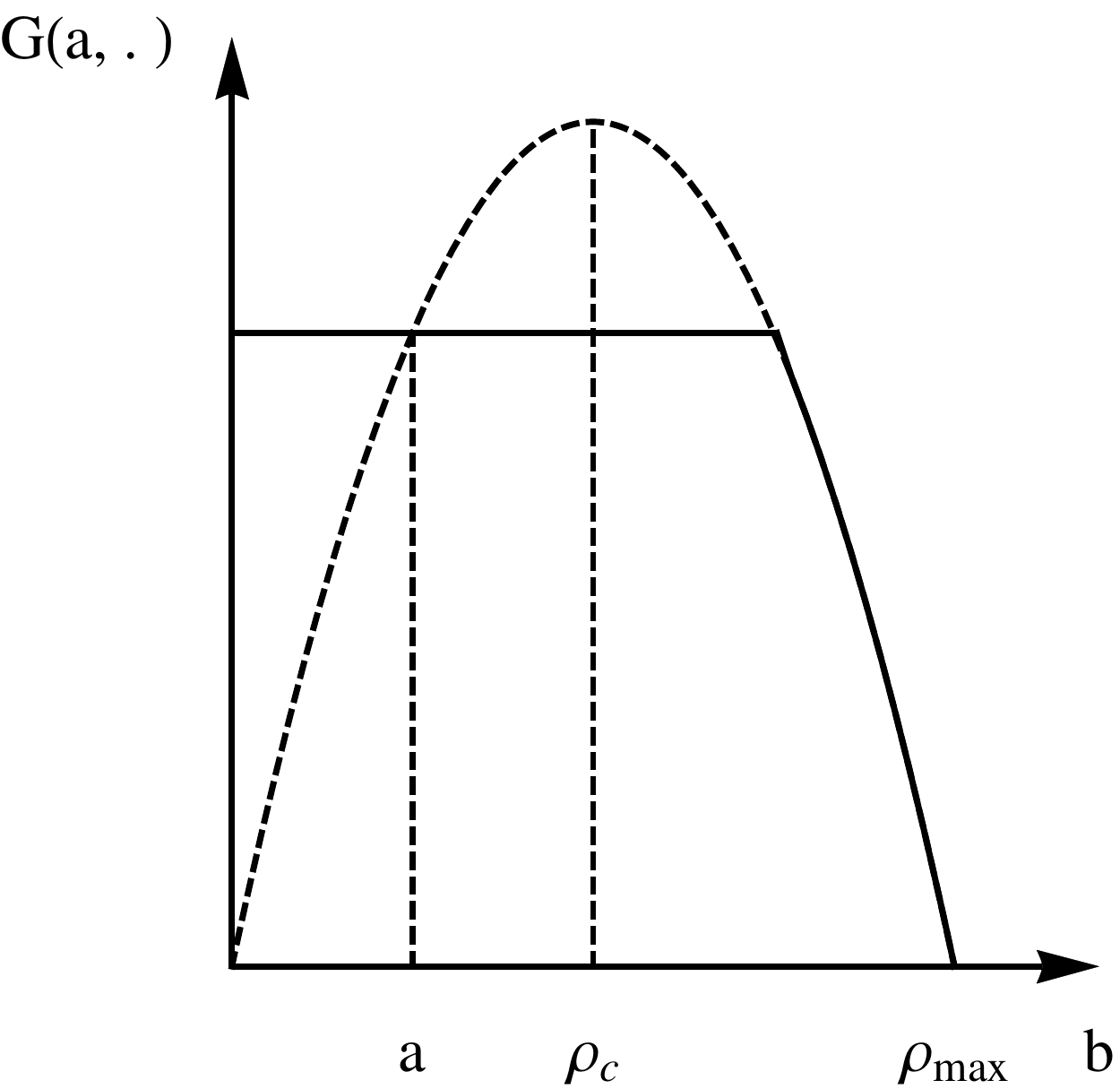}\quad
\includegraphics[width=.23\textwidth]{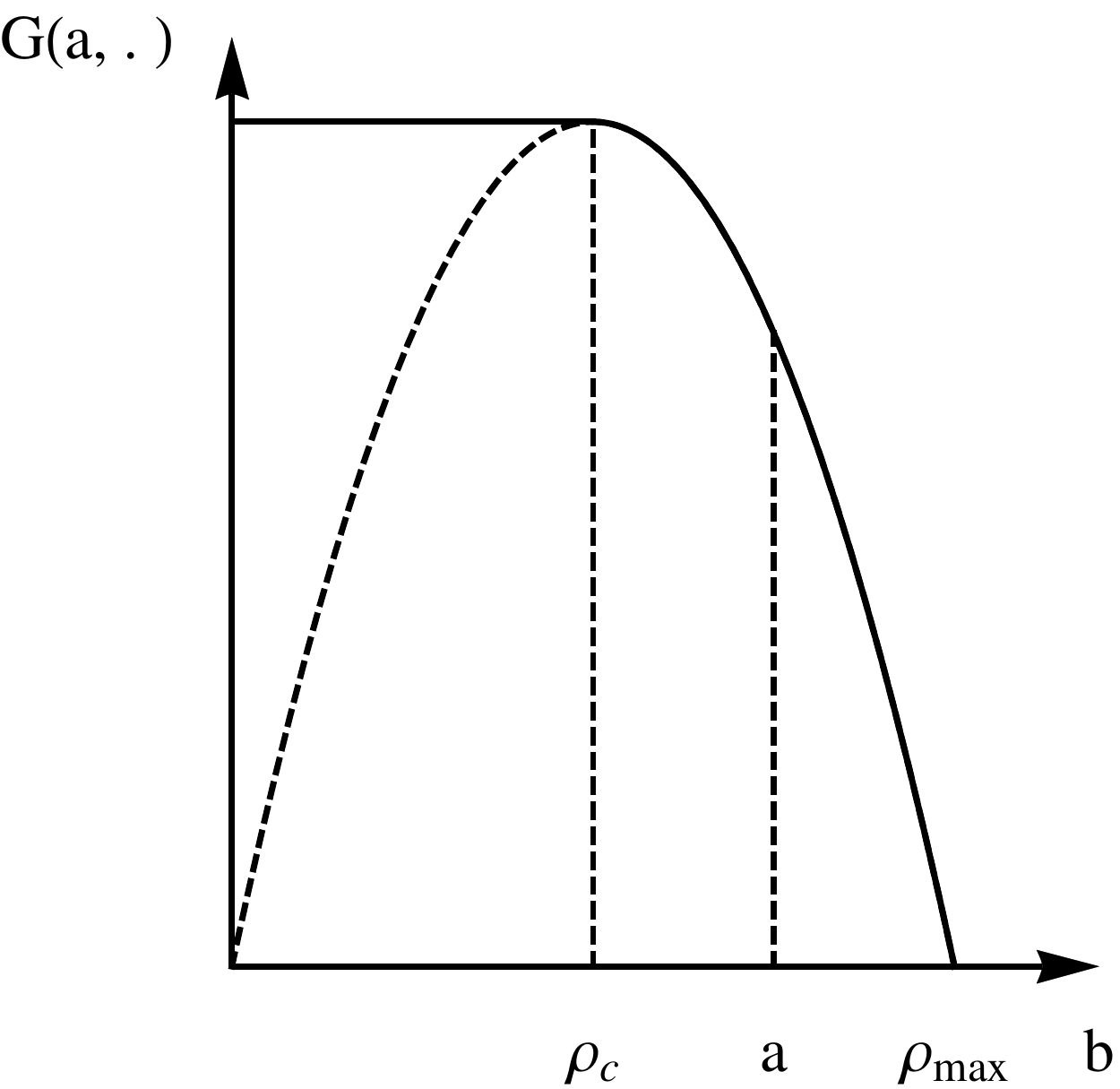}\quad
      \psfrag{G}[t,c]{$G(\cdot,b)$}
\includegraphics[width=.23\textwidth]{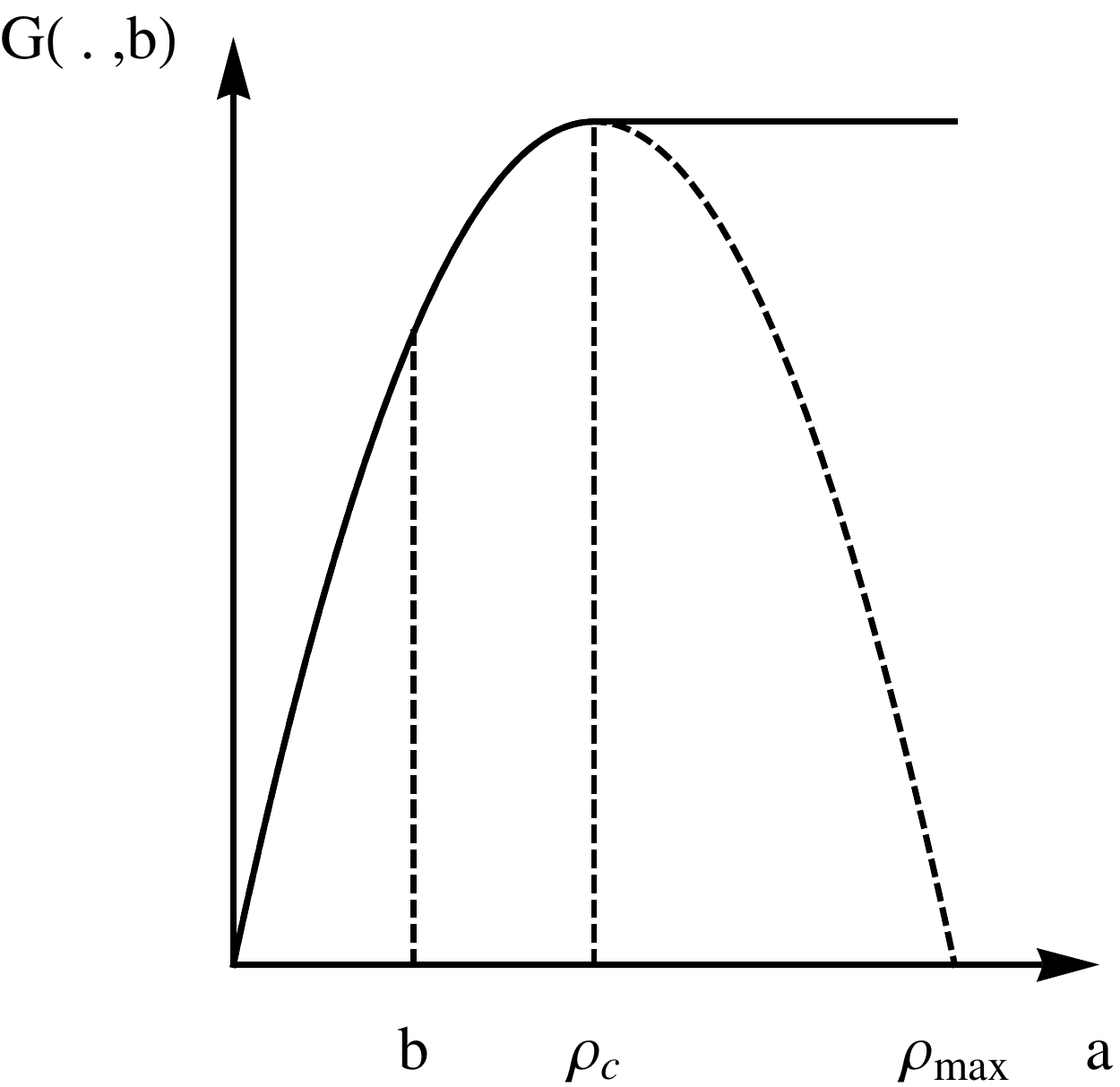}\quad
\includegraphics[width=.23\textwidth]{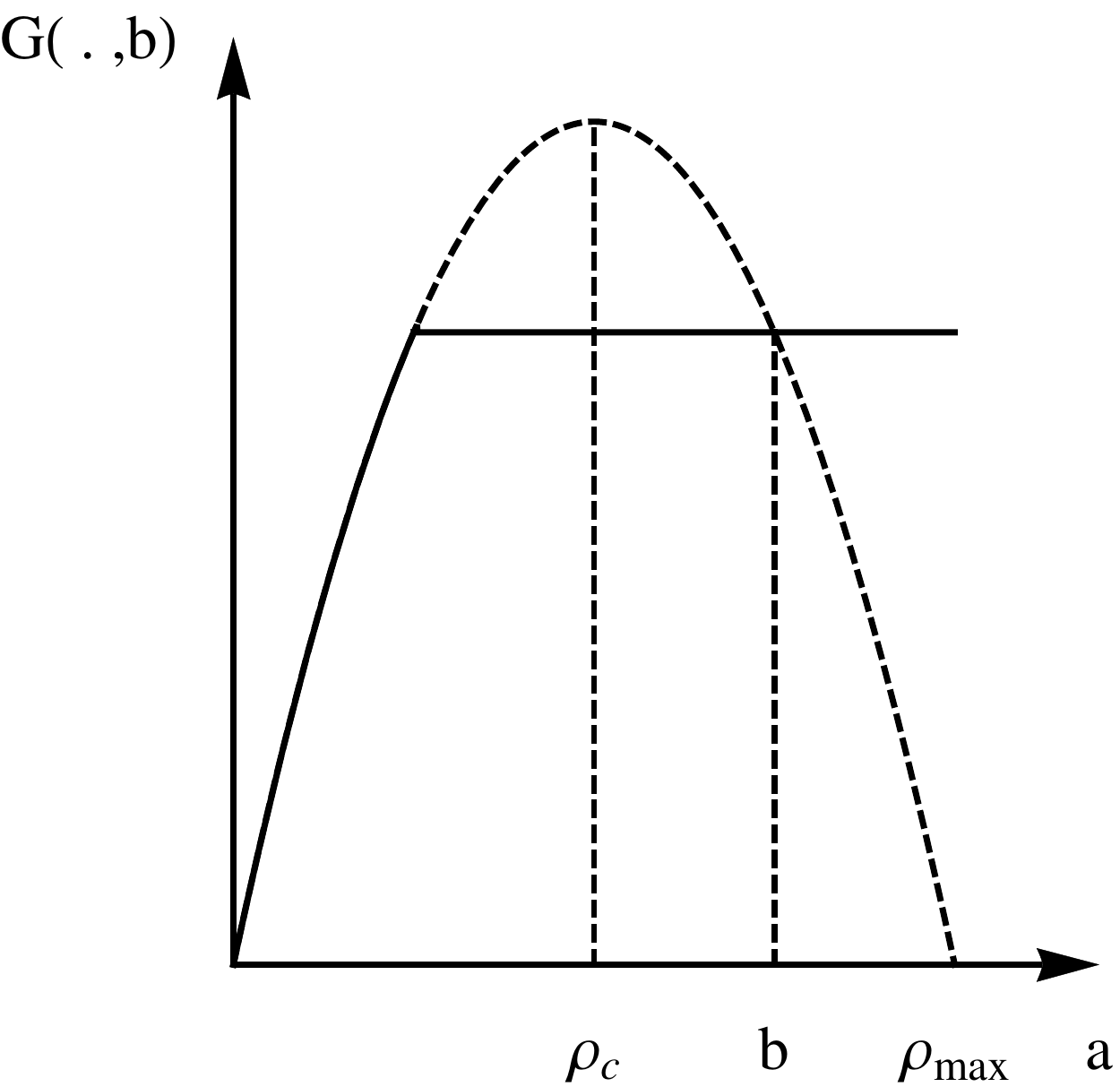}
\end{psfrags}}
\caption{The Godunov flux.}
\label{fig:G}
\end{figure}

\begin{definition}
The \emph{Godunov flux} related to a flux $f$ satisfying {\bf (F)} is the function which associates to any couple $(a,b) \in [0,\rho_{\max}]^2$ the value $f(\rho(t,0^-))=f(\rho(t,0^+))$ (the  equality holds due to the Rankine-Hugoniot condition), where $\rho$ is the Kruzhkov \cite{Kruzkov} entropy solution to the Riemann problem 
\begin{equation*}
 \begin{cases}
 \pt \rho + \px f(\rho) =0, &t>0,\ x\in\R,\\
 \rho(0,x) =
 \begin{cases}
  a& \text{if } x<0, \\
  b&\text{if } x>0,
 \end{cases}
 &x\in\R,
 \end{cases}
\end{equation*}
see \figurename~\ref{fig:G}.
The analytical expression of the Godunov flux is given by
\begin{equation*}\label{godunovflux}
G(a,b) =
\begin{cases}
 \min_{s\in[a,b]} \{f(s)\}&\text{if }a\leq b,\\
 \max_{s\in[b,a]} \{f(s)\}&\text{if }a\geq b.
\end{cases}
\end{equation*}
\end{definition}

We denote by $G_h$ the Godunov flux associated with the flux $f_h$, $h\in\mathsf{H}$.

Consider the initial boundary value problem
\begin{equation*}
\begin{cases}
\pt \rho + \px f(\rho) = 0,\qquad\text{for }(t,x)\in\R_+\times\R_-,\\
\rho(t,0)=\rho_b(t),\\
\rho(0,x)=\rho_0(x),
\end{cases}
\end{equation*}
and assume $\rho$ is a Kruzkov entropy solution in the interior of the half plane $\R_+\times\R_-$. Then $\rho$ satisfies the boundary condition in the sense of Bardos-LeRoux-N\'ed\'elec (see~\cite{Bardos}) if and only if $f(\rho(t,0-)) = G(\rho(t,0-),\rho_b(t))$.

Fix an initial condition $\vec\rho_0 = \left(\rho_{1,0},\cdots,\rho_{m+n,0}\right)\in\Lp\infty(\R_+\times\Gamma;[0,\rho_{\max}]^{m+n})$. We look for a function $\vec\rho = \left(\rho_1,\cdots,\rho_{m+n}\right)$ such that for every $h\in\mathsf{H}$, $\rho_h\in\Lp\infty(\R_+\times\Omega_h,[0,\rho_{\max}])$ is a weak entropy solution of
\begin{equation}
\label{eq:IBVP}
\begin{cases}
\pt\rho_h + \px f_h(\rho_h) = 0, &\quad\text{on }]0,T[\times\Omega_h,\\
\rho_h(t,0) = v_h(t), &\quad\text{on }]0, T[,\\
\rho_h(0,x) = \rho_{h,0}(x), &\quad\text{on }\Omega_h,
\end{cases}
\end{equation}
where $\vec v:\R_+\to[0,\rho_{\max}]^{m+n}$ is to be fixed in the sequel in order to ensure the conservation at the junction.

The condition~\eqref{eq:cons} is equivalent to ask for the traces $\rho_h(t,0\pm)$ to satisfy the boundary condition in the sense of Bardos-LeRoux-N\'ed\'elec
\begin{gather}
\label{BNL1}
f_h\left(\rho_h(t,0-)\right) = G_h\left(\rho_h(t,0-), v_h(t)\right),\qquad\text{if }h\in\mathsf{I};\\
\label{BNL2}
f_h\left(\rho_h(t,0+)\right) = G_h\left(v_h(t),\rho_h(t,0+)\right),\qquad\text{if }h\in\mathsf{J}.
\end{gather}

In order to describe the solution of~\eqref{eq:basic} we postulate (see~\cite{AC, AM})
\begin{equation}
\label{eq:p}
v_h(t) = p(t),\qquad\text{for a.e. }t\in\R_+,\quad\text{for all }h\in\mathsf{H}.
\end{equation}

The criterion for the choice of $p$ is equivalent to the condition~\eqref{eq:cons}, indeed, due to~\eqref{BNL1} and~\eqref{BNL2}, we can express~\eqref{eq:cons} in the following way
\begin{equation}
\label{eq:cons-godunov}
\sum_{i\in\mathsf{I}} G_i\left(\rho_i(t,0^-),p(t)\right)  = \sum_{j\in\mathsf{J}} G_j\left(p(t),\rho_j(t,0^+)\right),\qquad\text{for a.e. $t>0$}.
\end{equation}

We can now give the definition of admissible solution.
\begin{definition}
Given an initial condition $\vec\rho_0\in\Lp\infty(\Gamma;[0,\rho_{\max}]^{m+n})$, we say that $\vec\rho = \left(\rho_1,\cdots,\rho_{m+n}\right)$ in $\Lp\infty(\R_+\times\Gamma;[0,\rho_{\max}]^{m+n})$ is an {\em admissible solution} for the Cauchy problem at the junction~\eqref{eq:basic} associated with $\vec\rho_0$, if there exists a function $p\in\Lp\infty(\R_+:[0,\rho_{\max}])$ such that, for any $h\in\mathsf{H}$, $\rho_h$ is a weak solution of~\eqref{eq:IBVP} with $v_h$, $h\in\mathsf{H}$ chosen to fulfill~\eqref{eq:p}, and such that $\vec\rho$, $p$ fulfill~\eqref{eq:cons}.
\end{definition}

The authors provide a characterization of vanishing viscosity limits for the problem~\eqref{eq:basic} in terms of $m+n$ Dirichlet problems on $\Omega_h$, $h\in\mathsf{H}$ coupled by a transmission condition at the junction. 
To this aim, they introduce the {\em vanishing viscosity germ} (see~\cite{AKRvv, AKRL1, AC}), which can be identified by the set of all possible stationary admissible solutions to~\eqref{eq:basic} that are constant on each road of $\Gamma$.

The main result of~\cite{ACD} is summarized in the following theorem.

\begin{theorem}[Theorem 3.1 in~\cite{ACD}]
For any given initial condition $\vec\rho_0 = \left(\rho_{0,1}, \ldots, \rho_{0,m+n}\right)$ in $\Lp\infty(\Gamma; \R^{m+n})$ the problem \eqref{eq:basic} admits a unique admissible solution  $\vec\rho$ in $\Lp\infty(\R_+\times\Gamma ; [0,\rho_{\max}]^{m+n})$.

Moreover, 
if $\vec\rho$ and $\vec\rho'$ are two admissible solutions corresponding respectively to the initial condition $\vec\rho_0$ and $\vec\rho_0'$, then for all $M>0$ and $t<M/L$, where $L=\max \left\{ \|f_h'\|_{\Lp\infty([0,\rho_{\max}];\R)}\,|\, h=1,\ldots,m+n\right\}$,
 \begin{equation*}
  \begin{aligned}
   & \sum_{i=1}^{m} \int_{-(M-Lt)}^0 |\rho_{i}(t,x) - \rho_{i}'(t,x)|\,dx + \sum_{j=m+1}^{m+n} \int_0^{M-Lt} |\rho_{j}(t,x) - \rho_{j}'(t,x)|\,dx \\
        & \qquad\qquad \leq   \sum_{i=1}^{m} \int_{-M}^0 |\rho_{i,0}(x) - \rho_{i,0}'(x)|\,dx + \sum_{j=m+1}^{m+n} \int_0^M |\rho_{j,0}(x) - \rho_{j,0}'(x)|\,dx. \label{eq:L1loc-contraction}
  \end{aligned}
\end{equation*}
In particular, the map that associates to $\vec\rho_0$ the unique corresponding admissible profile $\vec{\rho}(t)$, is non-expansive with respect to the $\Lp1$ distance for all $t>0$.
\end{theorem}


\section{An explicit admissible solution at a merge}\label{sec:explicit}
In this section we compute an explicit solution for the problem consisting of two incoming and one outgoing arcs. We consider 
\[f(\rho)\equiv f_h(\rho) = \rho(1-\rho)\]
as the flux for each arc. As initial condition, we choose 
\[\rho_{1,0}(x)=\chi_{[-1/2, 0]}(x),\quad \rho_{2,0}(x)=3/4\chi_{[-1/4,0]}(x),\quad \rho_{3,0}(x)=0.\]

The exact solution is obtained by an explicit analysis of the wave-front interactions, with computer assisted computation of the front slopes and interaction times. 

At time $t=0$, let $p_1\approx 0.85$ be the solution of
\[G_1(1,p_1)+G_2(3/4,p_1)=G_3(p_1,0),\]
then, on $\Omega_1$ a rarefaction $\mathfrak{R}_{O,1}$ starts from $O(0,0)$ and its values are given by
\[\mathfrak{R}_{O,1}=\frac{1}{2}\left(1-\frac{x}{t}\right),\quad\text{for }-t\le x \le -\frac{\sqrt{2}}{2}t.\]
On $\Omega_2$ starts the backward shock $\mathfrak{S}_{O,2}$ given by
\[\mathfrak{S}_{O,2}:\,\dot{x}(t) = \sigma\left(\frac{3}{4}, p_1\right),\quad x(0)=0.\]
On $\Omega_3$ a rarefaction starts from $O(0,0)$ and its values are given by
\[\mathfrak{R}_{O,3} = \frac{1}{2}\left(1-\frac{x}{t}\right),\quad\text{for }0\le x \le t.\]
On $\Omega_2$, let $C(x_C,t_C)$ be the point where the shock $\mathfrak{S}_{B,2}:\,x(t)=-\frac{1}{4} + \frac{t}{4}$ originated from $B(-1/4,0)$ interacts with the shock $\mathfrak{S}_{O,2}$. As a result, from $C$ starts a shock given by
\[ \mathfrak{S}_{C,2}:\,\dot{x}(t)=\sigma\left(0,\,p_1\right),\quad x(t_C)=x_C,
\]
and reaches the junction $x=0$ at time $t=t_F=3/2$ that corresponds to the time at which the second incoming arc becomes empty.
On $\Omega_1$, in $D(-1/2,1/2)$, the stationary shock $\mathfrak{S}_{A,1}$ originated from $A(-1/2,0)$ interacts with the rarefaction $\mathfrak{R}_{O,1}$. As a result, from $D$ starts a shock $\mathfrak{S}_{D,1}$ given by
\begin{equation*}
\mathfrak{S}_{D,1}:\, \dot{x}(t)=\sigma\left(0,\,\mathfrak{R}_{O,1}\left(t,x(t)\right)\right)
,\quad x(1/2)=-1/2.
\end{equation*}
Let $E(x_E,t_E)$ be the intersection between $\mathfrak{S}_{D,1}$ and $x(t)=-(\sqrt{2}/2)\,t$. From this point starts a forward shock
\[
\mathfrak{S}_{E,1}: \, \dot{x}(t)=\sigma\left(0,\,p_1\right)
,\quad x(t_E)=x_E.
\]
Let $p_2=1/2$ be the solution of
\[G_1\left(p_1,p_2\right)+G_2\left(0,p_2\right)=G_3\left(p_2,\frac{1}{2}\right),\quad\text{for }t > t_F.\]
Therefore, a rarefaction appears on $\Omega_1$:
\begin{equation*}
\mathfrak{R}_{F,1}(t,x)=\frac{1}{2}\left(1-\frac{x}{t-t_F}\right),\qquad\text{for}\quad -\frac{\sqrt{2}}{2}(t-t_F)< x \le 0.
\end{equation*}Let $G$ be the point where $\mathfrak{S}_{E,1}$ and $\mathfrak{R}_{F,1}$ interact. From this point starts a forward shock $\mathfrak{S}_{G,1}$, with left state $\rho = 0$, which reaches the junction at time $t_H\approx 2.75$, then $\Omega_1$ is empty.
Finally, on $\Omega_3$ at time $t_H$ starts a shock which interacts with the rarefaction $\mathfrak{R}_{O,3}$ generating the additional shock 
\[
\mathfrak{S}_{H,3}: \, \dot{x}(t)=\sigma\left(0,\,\mathfrak{R}_{O,3}\left(t,x(t)\right)\right)
,\quad x(t_H)=0.
\]

\section{Finite volumes numerical scheme}~\label{sec:volume-scheme}

We fix a constant space step $\Delta x$. For $\ell \in \Z$ and $h \in \mathsf{H}$, we set $x_{\ell}^h = \ell \Delta x$. We define the cell centers $x_{\ell+\frac{1}{2}}^h = (\ell +\frac{1}{2})\Delta x$ for $\ell\in\Z$ and consider the uniform spatial mesh on each $\Omega_h$
\begin{equation}
\label{eq:spatial-mesh}
\bigcup_{\ell \le -1}{(x_{\ell}^i, x_{\ell+1}^i)},\quad i\in\mathsf{I},\qquad
\bigcup_{\ell \ge 0}{(x_{\ell}^j, x_{\ell+1}^j)},\quad j\in\mathsf{J},
\end{equation}
so that the position of the junction $x = 0$ corresponds to $x_0^h$ for each edge. Then we fix a constant time step $\Delta t$ satisfying the CFL condition
\begin{equation}
\label{eq:CFL}
\Delta t \max_{h}\{L_h\} \le \frac{\Delta x}{2},
\end{equation}
and for $s \in \N$ we define the time discretization $t^s=s \Delta t$. At each time $t^s$, $\rho_{\ell+\frac{1}{2}}^{h,s}$ represents an approximation of the main value of the solution on the interval $[x_{\ell}^h,x_{\ell+1}^h)$, $\ell\in\Z$, along the $h$-th arc. We initialize the scheme by discretizing the initial conditions
\begin{equation}
\label{eq:init-dicret}
\rho_{\ell + \frac{1}{2}}^{h,0} = \frac{1}{\Delta x} \int_{x_{\ell}^h}^{x_{\ell+1}^h}\rho_h^0(x)\,dx,
\end{equation}
for all $h \in\mathsf{H}$ and for $\ell \le -1$ if $h \in\mathsf{I}$, $\ell \ge 0$ if $h \in\mathsf{J}$.

For each $s\in\N$, at all cell interfaces $x_{\ell}^h$ with $\ell \ne 0$ we consider the standard Godunov flux $G_h$ corresponding to the flux $f_h$. At the junction $x_0^h$ we take on each arc $\Omega_h$ the Godunov flux corresponding to the admissible solution of the Riemann problem at the junction, defined as in \cite{ACD}, which we compute by means of a two-step procedure:
\begin{enumerate}
	\item[(i)] find
		\begin{equation}
		\label{eq:find-p}
		p^s \in [0,\rho_{\max}]\quad\text{s.t. }\sum_{i\in\mathsf{I}} G_i(\rho_{-\frac{1}{2}}^{i,s}, p^s) = \sum_{j\in\mathsf{J}} G_j(p^s, \rho_{\frac{1}{2}}^{j,s}),
		\end{equation}
	\item[(ii)] compute
		\begin{equation}
		\label{eq:update}
		\rho_{\ell+\frac{1}{2}}^{h,s+1} = \rho_{\ell+\frac{1}{2}}^{h,s} - \frac{\Delta t}{\Delta x} \left(\mathcal{F}_{\ell+1}^{h,s}-\mathcal{F}_{\ell}^{h,s}\right),
		\end{equation}
		where
		\begin{equation}
		\label{eq:num-flux}
		\mathcal{F}_{\ell}^{h,s}=
		\begin{cases}
		G_h\left(\rho_{\ell-\frac{1}{2}}^{h,s}, \rho_{\ell+\frac{1}{2}}^{h,s}\right),\quad&\text{if }h \in\mathsf{I}\text{ and }\ell \ge -1\text{ or }h \in\mathsf{J}\text{ and }\ell \ge 1,\\
		G_h(\rho_{-\frac{1}{2}}^{h,s}, p^s),\quad&\text{if }h \in\mathsf{I}\text{ and }\ell = 0,\\
		G_h(p^s, \rho_{\frac{1}{2}}^{h,s}),\quad&\text{if }h \in\mathsf{J}\text{ and }\ell=0.
		\end{cases}
		\end{equation}
\end{enumerate}

The choice of the Godunov's flux is motivated by the fact that all admissible stationary solutions are exact solutions for such scheme. However, one can use any other numerical flux that is monotone, consistent and Lipschitz continuous.

A convergence result for the scheme~\eqref{eq:find-p}-\eqref{eq:update}-\eqref{eq:num-flux} can be found in~\cite{ACD}.

\section{Validation of the numerical scheme}\label{sec:validation}
We propose here to validate the numerical scheme \eqref{eq:find-p}-\eqref{eq:update}-\eqref{eq:num-flux} making a comparison with the explicit solution computed in the Section~\ref{sec:explicit}; by comparison with the solution of Riemann problems and by showing the consistence of the scheme with respect to a network consisting of one single arc.

\begin{figure}[!htbp]%
\centering
	\begin{subfigure}[b]{.328\textwidth}
			\includegraphics[width=\textwidth]{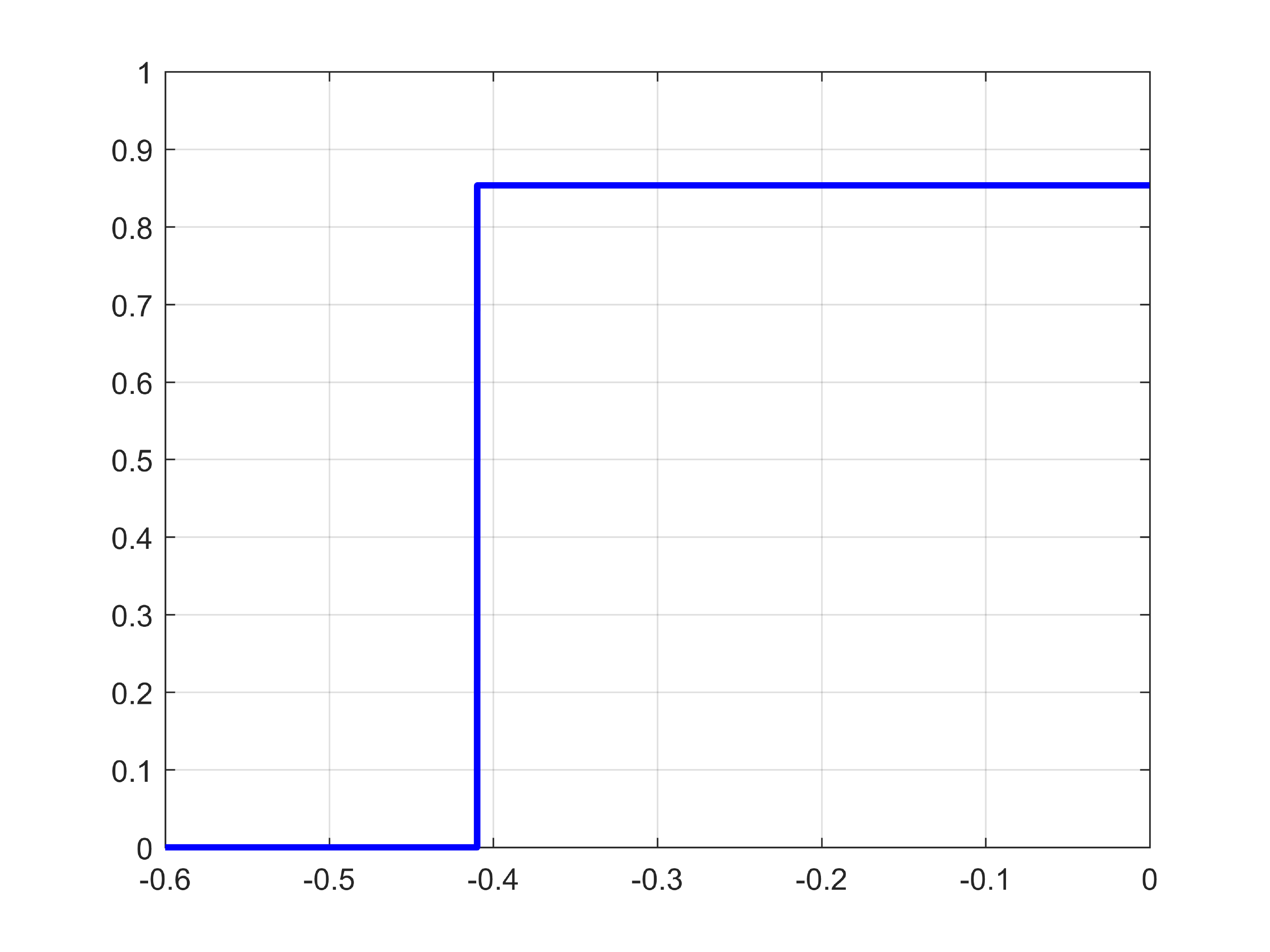}
			\caption*{$x\mapsto\rho_{1,\Delta}(1.2,x)$}
	\end{subfigure}
	\begin{subfigure}[b]{.328\textwidth}
			\includegraphics[width=\textwidth]{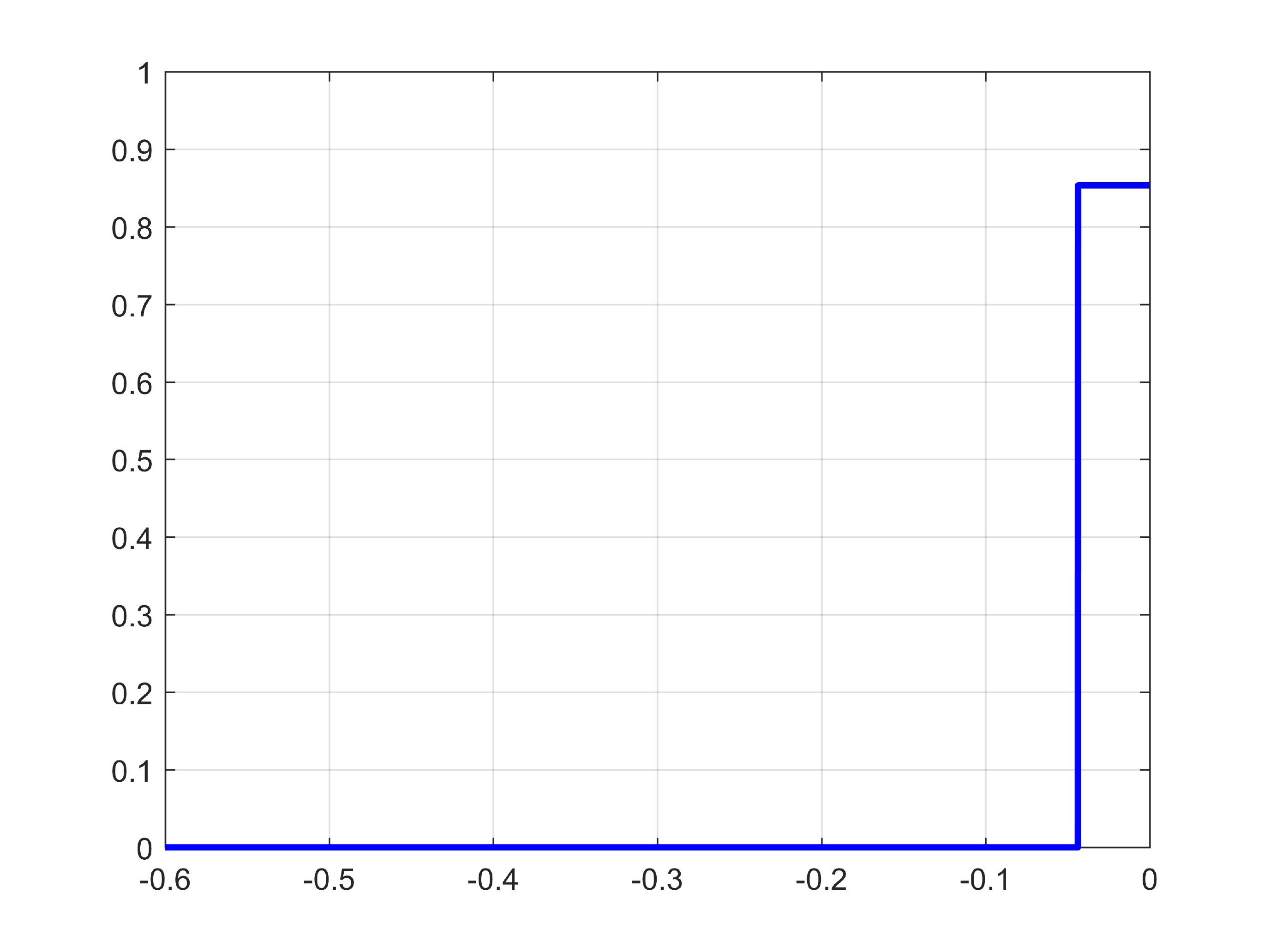}
			\caption*{$x\mapsto\rho_{2,\Delta}(1.2,x)$}
	\end{subfigure}
	\begin{subfigure}[b]{.328\textwidth}
			\includegraphics[width=\textwidth]{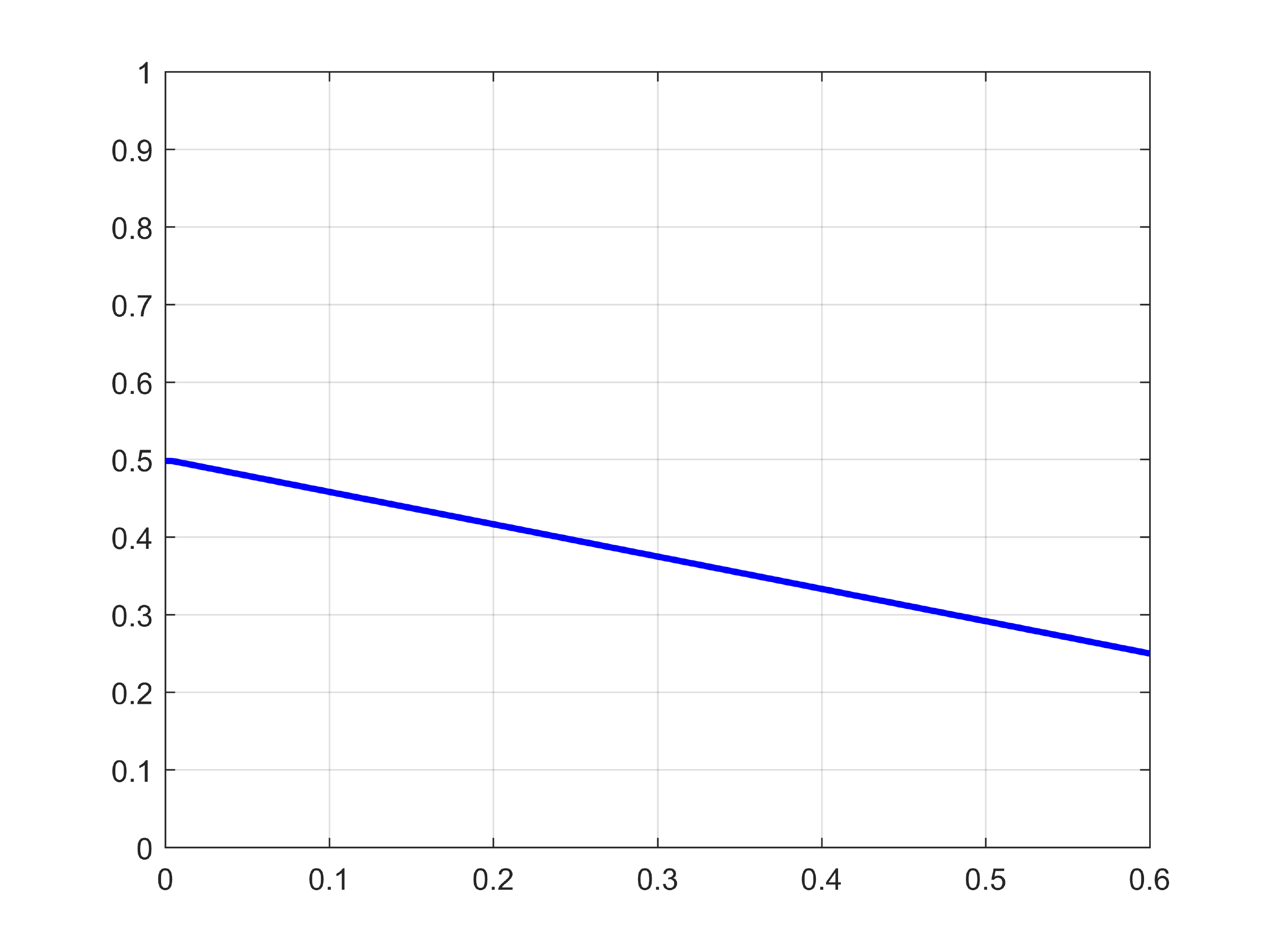}
			\caption*{$x\mapsto\rho_{3,\Delta}(1.2,x)$}
	\end{subfigure}
	\\
	\begin{subfigure}[b]{.328\textwidth}
			\includegraphics[width=\textwidth]{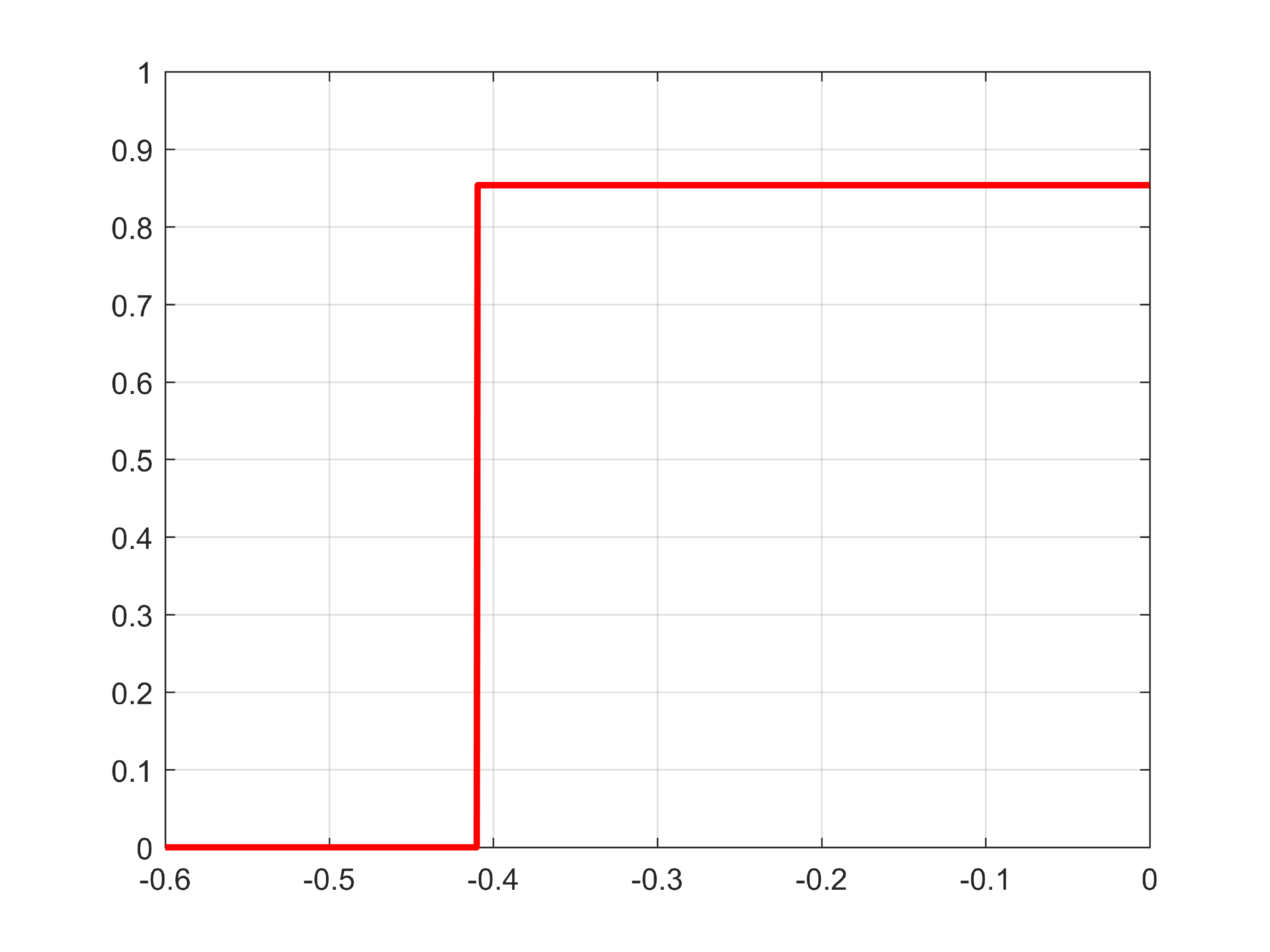}
			\caption*{$x\mapsto\rho_{1}(1.2,x)$}
	\end{subfigure}
	\begin{subfigure}[b]{.328\textwidth}
			\includegraphics[width=\textwidth]{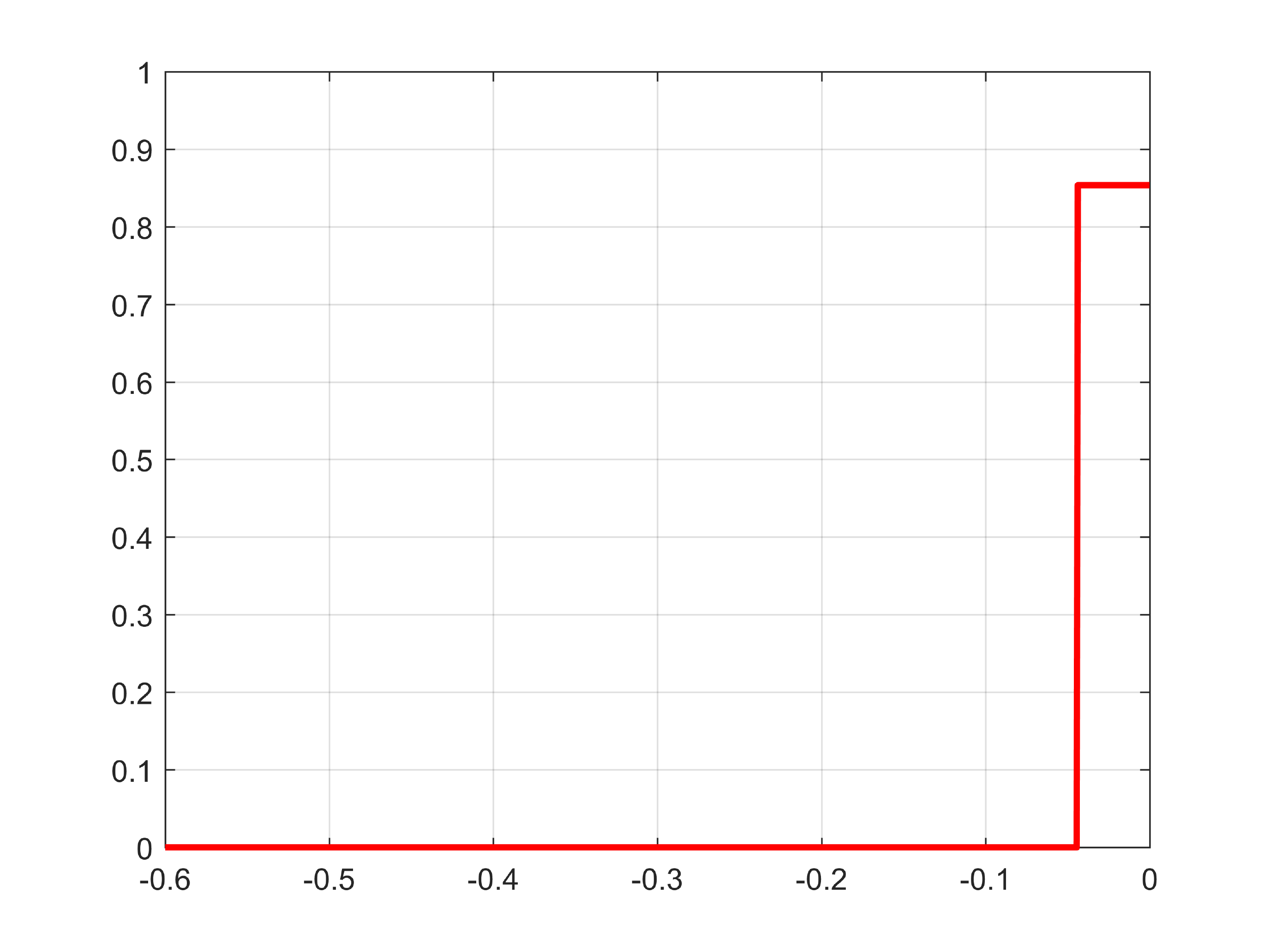}
			\caption*{$x\mapsto\rho_{2}(1.2,x)$}
	\end{subfigure}
	\begin{subfigure}[b]{.328\textwidth}
			\includegraphics[width=\textwidth]{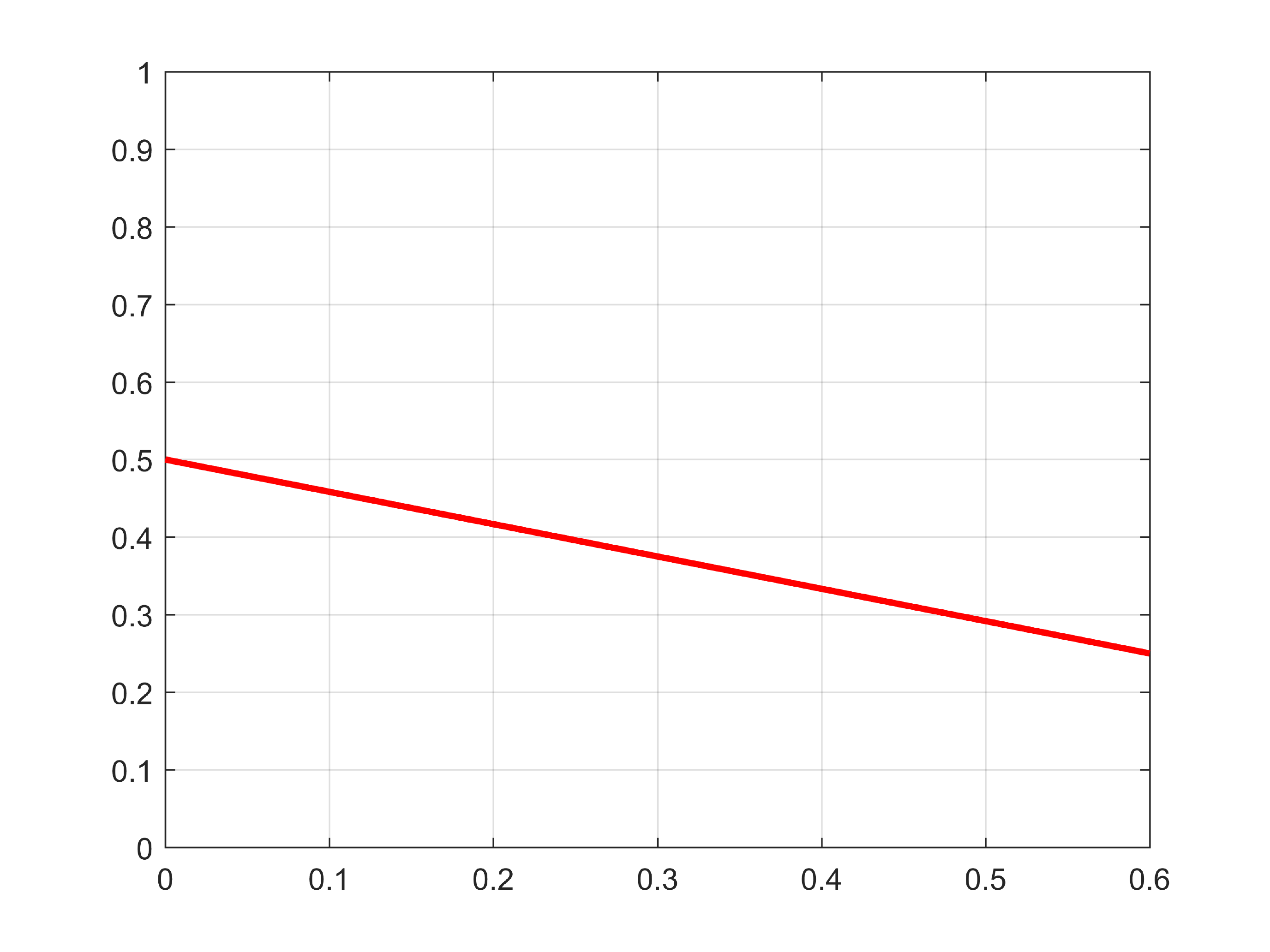}
			\caption*{$x\mapsto\rho_{3}(1.2,x)$}
	\end{subfigure}
	\\
	\begin{subfigure}[b]{.328\textwidth}
			\includegraphics[width=\textwidth]{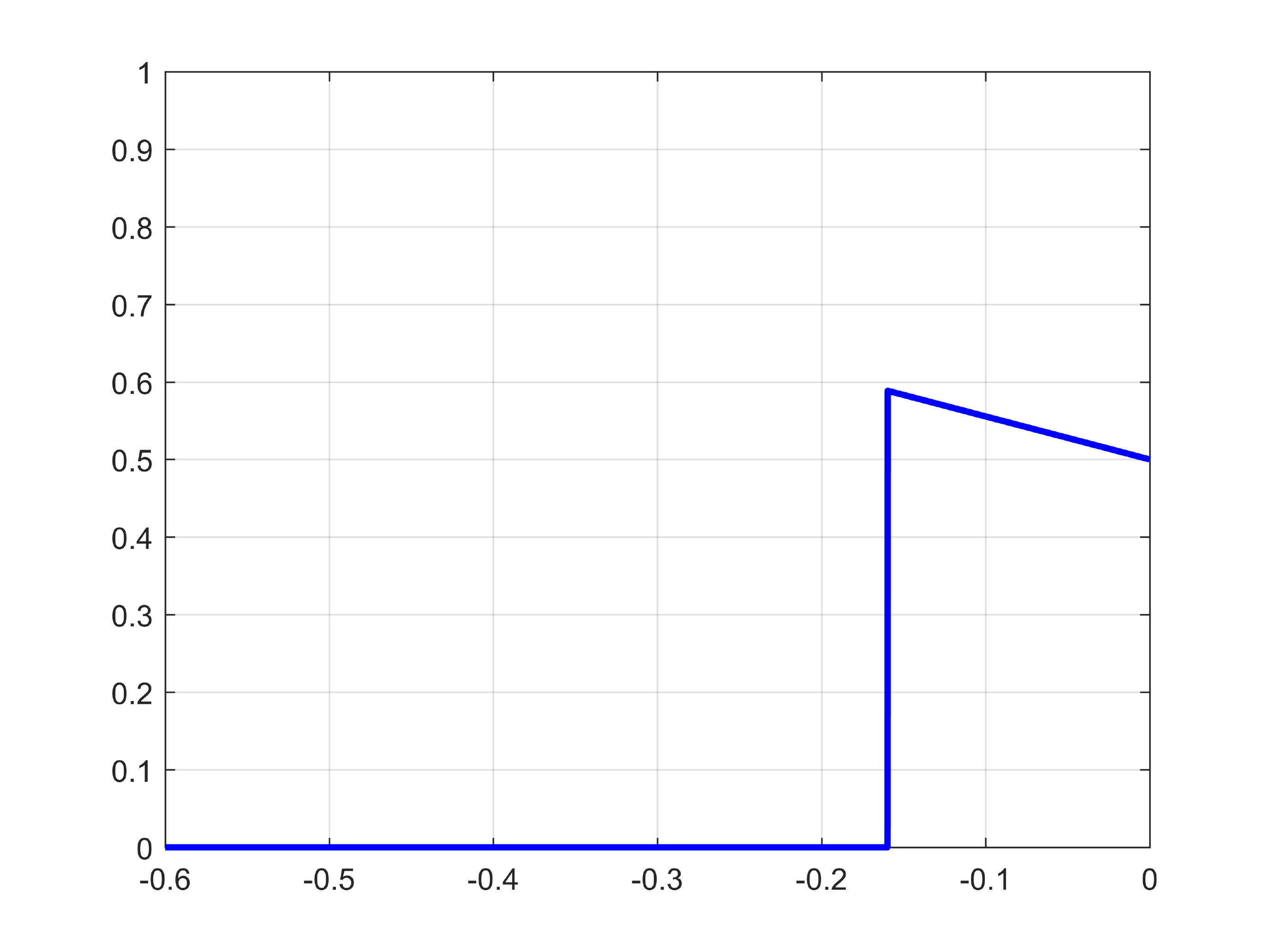}
			\caption*{$x\mapsto\rho_{1,\Delta}(2.4,x)$}
	\end{subfigure}
	\begin{subfigure}[b]{.328\textwidth}
			\includegraphics[width=\textwidth]{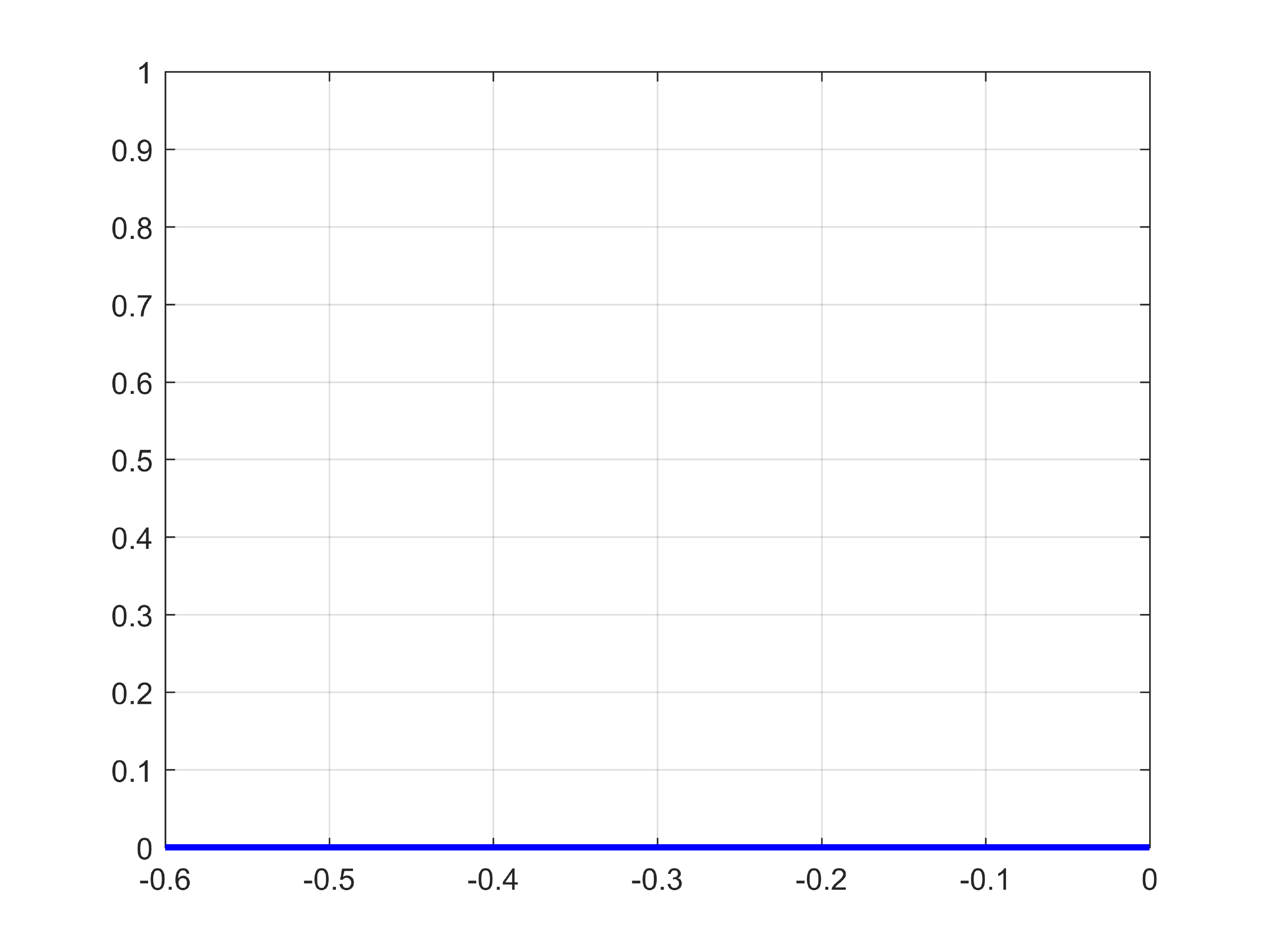}
			\caption*{$x\mapsto\rho_{2,\Delta}(2.4,x)$}
	\end{subfigure}
	\begin{subfigure}[b]{.328\textwidth}
			\includegraphics[width=\textwidth]{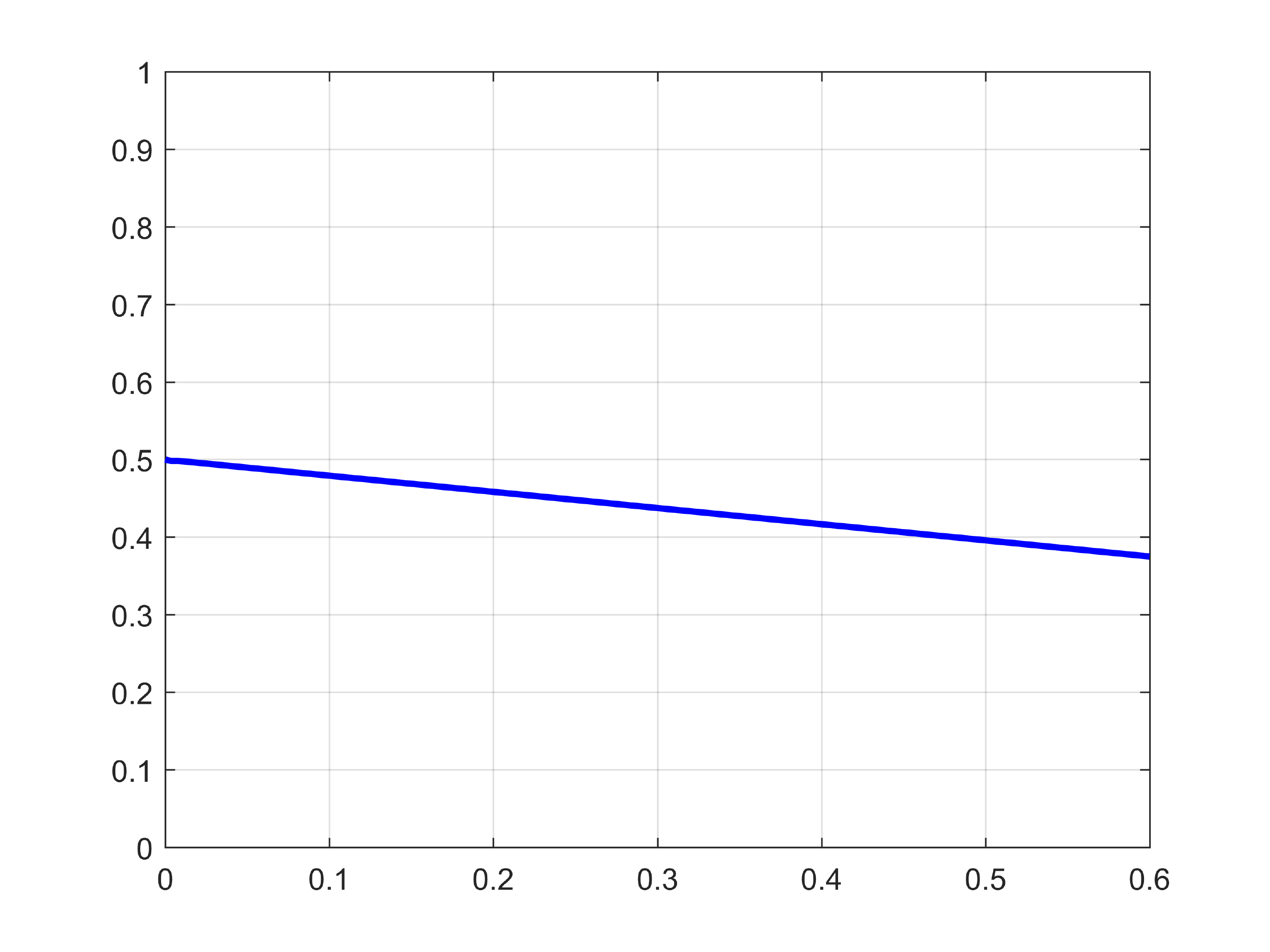}
			\caption*{$x\mapsto\rho_{3,\Delta}(2.4,x)$}
	\end{subfigure}
	\\
	\begin{subfigure}[b]{.328\textwidth}
			\includegraphics[width=\textwidth]{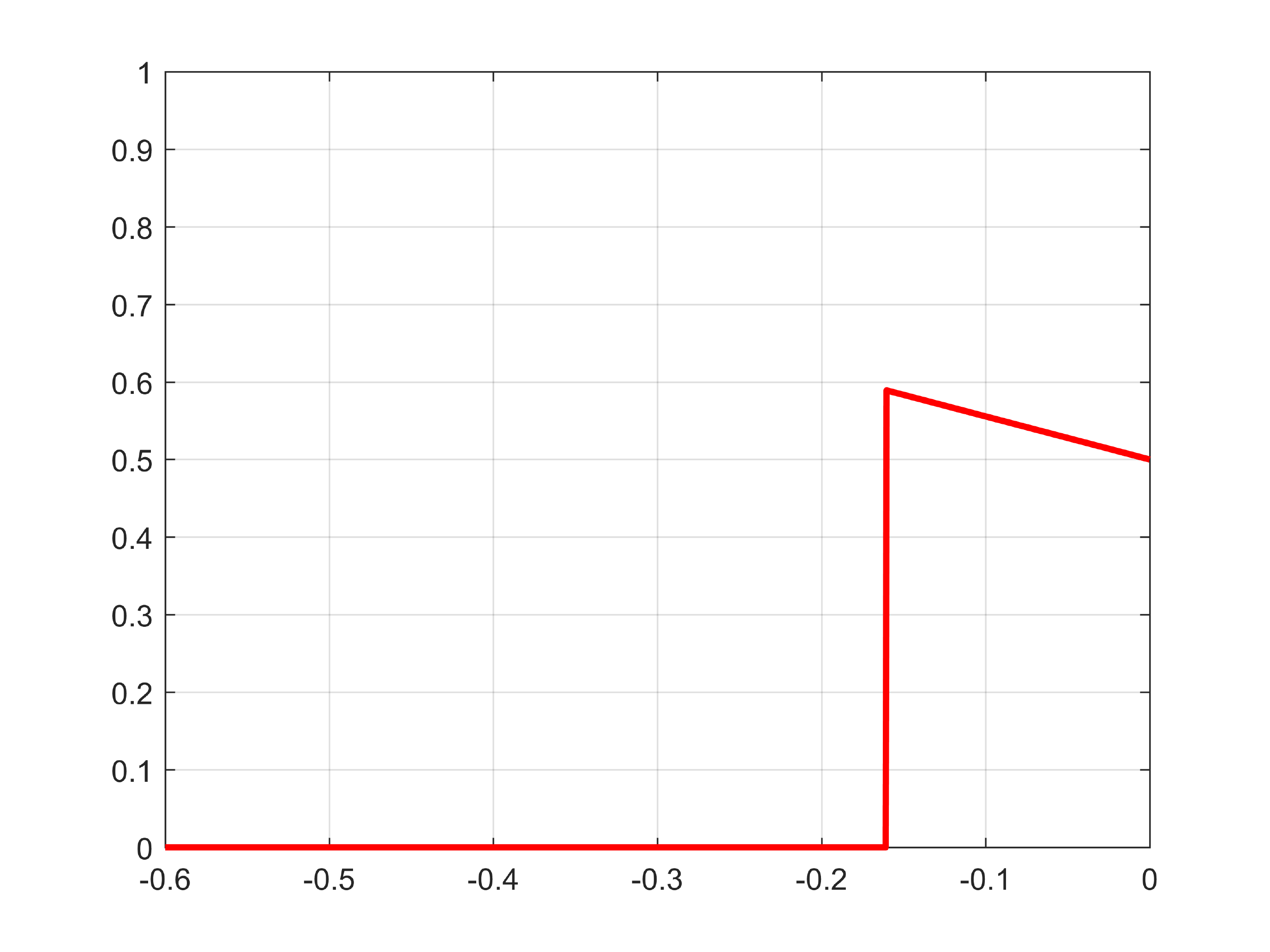}
			\caption*{$x\mapsto\rho_{1}(2.4,x)$}
	\end{subfigure}
	\begin{subfigure}[b]{.328\textwidth}
			\includegraphics[width=\textwidth]{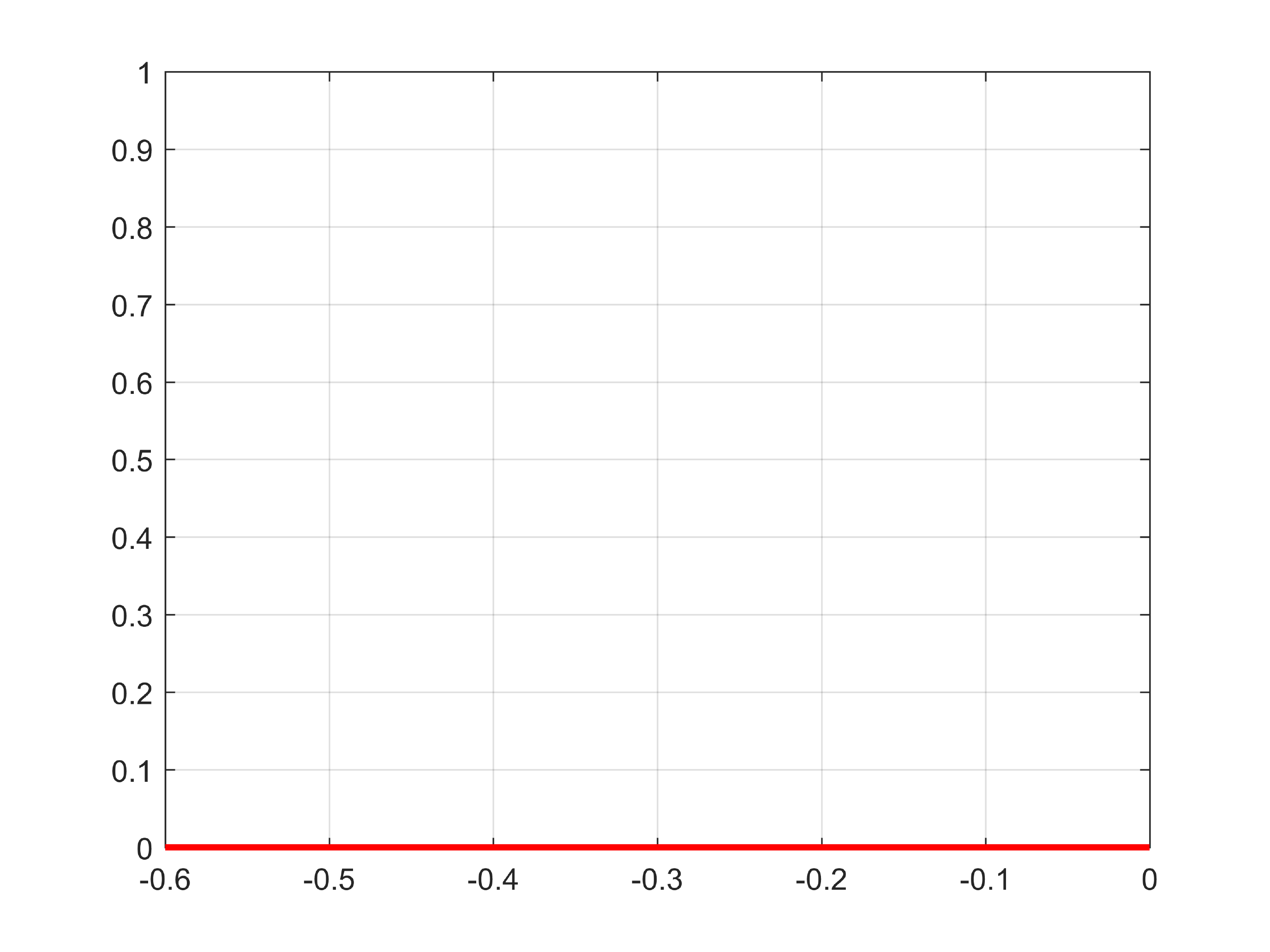}
			\caption*{$x\mapsto\rho_{2}(24,x)$}
	\end{subfigure}
	\begin{subfigure}[b]{.328\textwidth}
			\includegraphics[width=\textwidth]{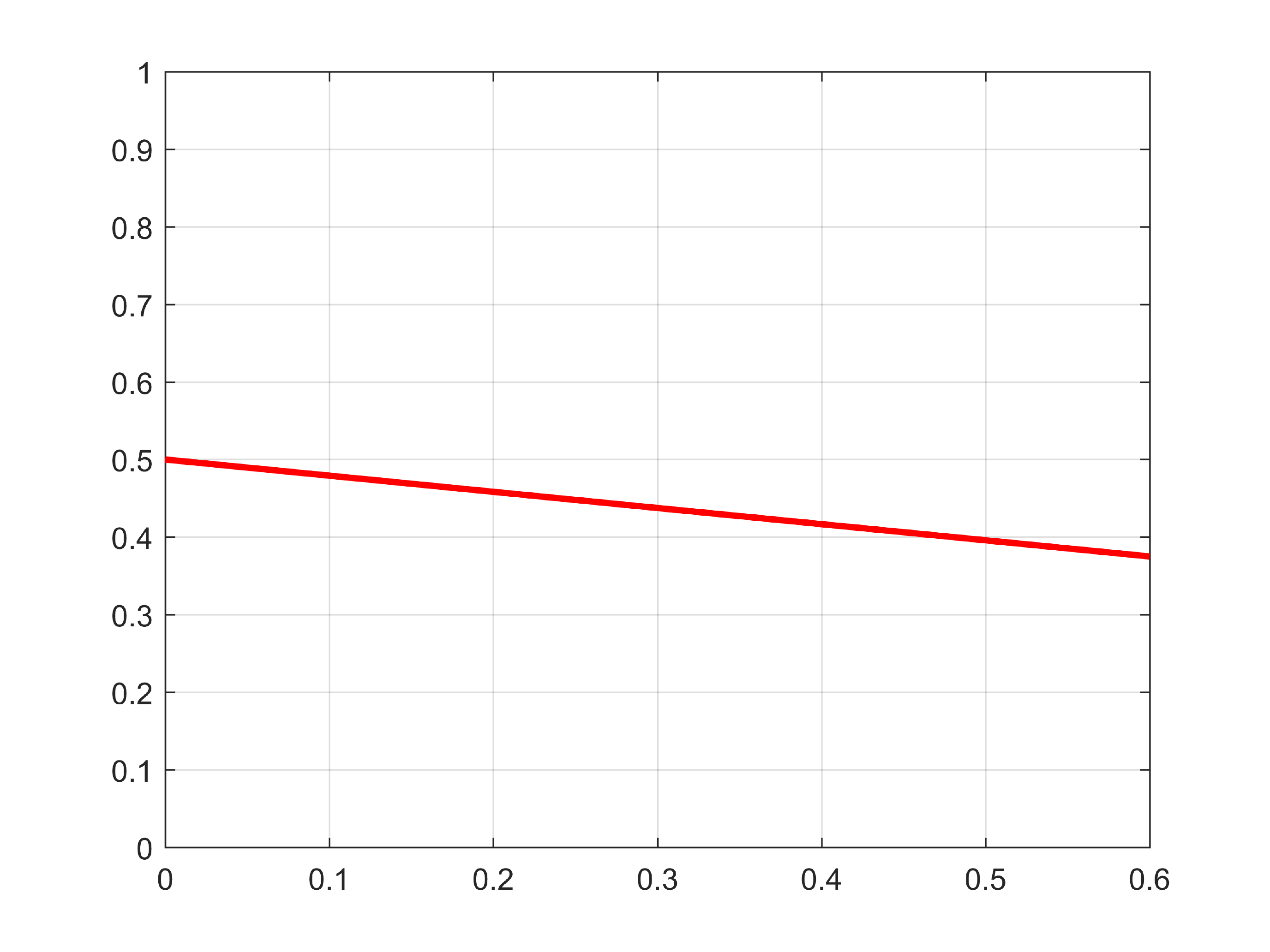}
			\caption*{$x\mapsto\rho_{3}(2.4,x)$}
	\end{subfigure}
\caption[The comparison between the explicit and numerical solution at times $t=1.2$ and $t=2.4$]{With reference to the simulation of Section~\ref{sec:validation}: comparison between the explicit solution $\vec\rho$ and the numerical one $\vec\rho_{\Delta}$ at times $t=1.2$ and $t=2.4$.}
\label{fig:explicit}
\end{figure}

We consider the explicit solution to~\eqref{eq:basic} constructed in Section~\ref{sec:explicit}. The setup for the simulation is as follows. We consider $[-3/5,0]$ as domain of computation for the incoming arcs and $[0, 3/5]$ for the outgoing one, and $\Delta x = 0.5\times 10^{-4}$, $\Delta t = 0.25 \times 10^{-4}$ as space and time step, respectively. A qualitative comparison between the numeric solution $x\mapsto\rho_{\Delta}(t,x)$ and the explicit solution $x\mapsto\rho(t,x)$ at different fixed times $t$ is shown in Figure~\ref{fig:explicit}.

Additionally, we perform a convergence analysis for this test. We introduce the relative $\Lp1$-error respectively for the whole network, for the incoming and for the outgoing arcs at a given time $t^s$ as follows
\begin{align*}
E_{\Lp1}^{s,\Gamma}&=\frac{\sum_{h\in\mathsf{H}}\sum_{\ell}|\rho_h(t^s,x_{\ell})-\rho_{\ell}^{h,s}|}{\sum_{h=1}^3\sum_{\ell}|\rho_h(t^s,x_{\ell})|},
\\
E_{\Lp1}^{s,\mathsf{I}}&=\frac{\sum_{i\in\mathsf{I}}\sum_{\ell}|\rho_i(t^s,x_{\ell})-\rho_{\ell}^{i,s}|}{\sum_{i=1}^2\sum_{\ell}|\rho_i(t^s,x_{\ell})|},
\\
E_{\Lp1}^{s,\mathsf{J}}&=\frac{\sum_{j\in\mathsf{J}}\sum_{\ell}|\rho_j(t^s,x_{\ell})-\rho_{\ell}^{j,s}|}{\sum_{j\in\mathsf{J}}\sum_{\ell}|\rho_j(t^s,x_{\ell})|}.
\end{align*}
Table~\ref{tab:rateconv} depicts the relative $\Lp1$-error with respect to the space step at the fixed time $t=2.4$. The time step is fixed to $\Delta t = 0.25\times 10^{-4}$. Since we are dealing with a first order scheme approximating discontinuous solutions, the sub-linear convergence rate found results expected.

\begin{table}%
\centering%
\renewcommand\arraystretch{1.3}
\begin{tabular}{|c|c|c|c|c|c|c|c|
}
\hline
Number & 
& Rate of &
&Rate of&
&Rate of
\\
of cells&  $E_{\Lp1}^{s,\Gamma}$ &convergence& $E_{\Lp1}^{s,\mathsf{I}}$ &convergence& $E_{\Lp1}^{s,\mathsf{J}}$&convergence
\\
per arc&&&&&&
\\
\hline
$60$&$6.5374 \times 10^{-2}$&-&$2.1928\times10^{-1}$&-&$1.4155\times10^{-2}$&-
\\
\hline
$120$&$3.4281\times 10^{-2}$&$0.9313$&$1.1380\times10^{-1}$&$0.9464$&$7.8554\times10^{-3}$&$0.8496$
\\
\hline
$600$&$7.6754\times 10^{-3}$&$0.9302$&$2.4933\times10^{-2}$&$0.9441$&$1.9468\times10^{-3}$&$0.8625$
\\
\hline
$1200$&$4.8890\times10^{-3}$&$0.8800$&$1.6393\times10^{-2}$&$0.8831$&$1.0579\times10^{-3}$&$0.8660$
\\
\hline
$6000$&$1.9875\times 10^{-3}$&$0.7721$&$7.1933\times 10^{-3}$&$0.7575$&$2.5294\times10^{-4}$&$0.8738$
\\
\hline
$12000$&$1.6804\times 10^{-3}$&$0.7034$&$6.3143\times 10^{-3}$&$0.6815$&$1.3587\times10^{-4}$&$0.8774$
\\
\hline
\end{tabular}
\renewcommand\arraystretch{1}
\caption{Relative $\Lp1$-error at time $t=2.4$ computed in Section~\ref{sec:validation}.}
\label{tab:rateconv}
\end{table}

\subsection{Riemann problem for a $2$-$1$ merge.}\label{sec:2-1}
We consider a network consisting of three edges and one junction, with two incoming and one outgoing arcs. We consider $[-\frac{1}{2}, 0]$ as domain of computation for the incoming arcs and $[0,\frac{1}{2}]$ for the outgoing one, and we take a normalized flux $f(\rho) = \rho \left(1 - \rho \right)$ for each arc. 

\begin{figure}
	\centering
	\begin{subfigure}[b]{0.328\textwidth}
			\includegraphics[width=\textwidth]{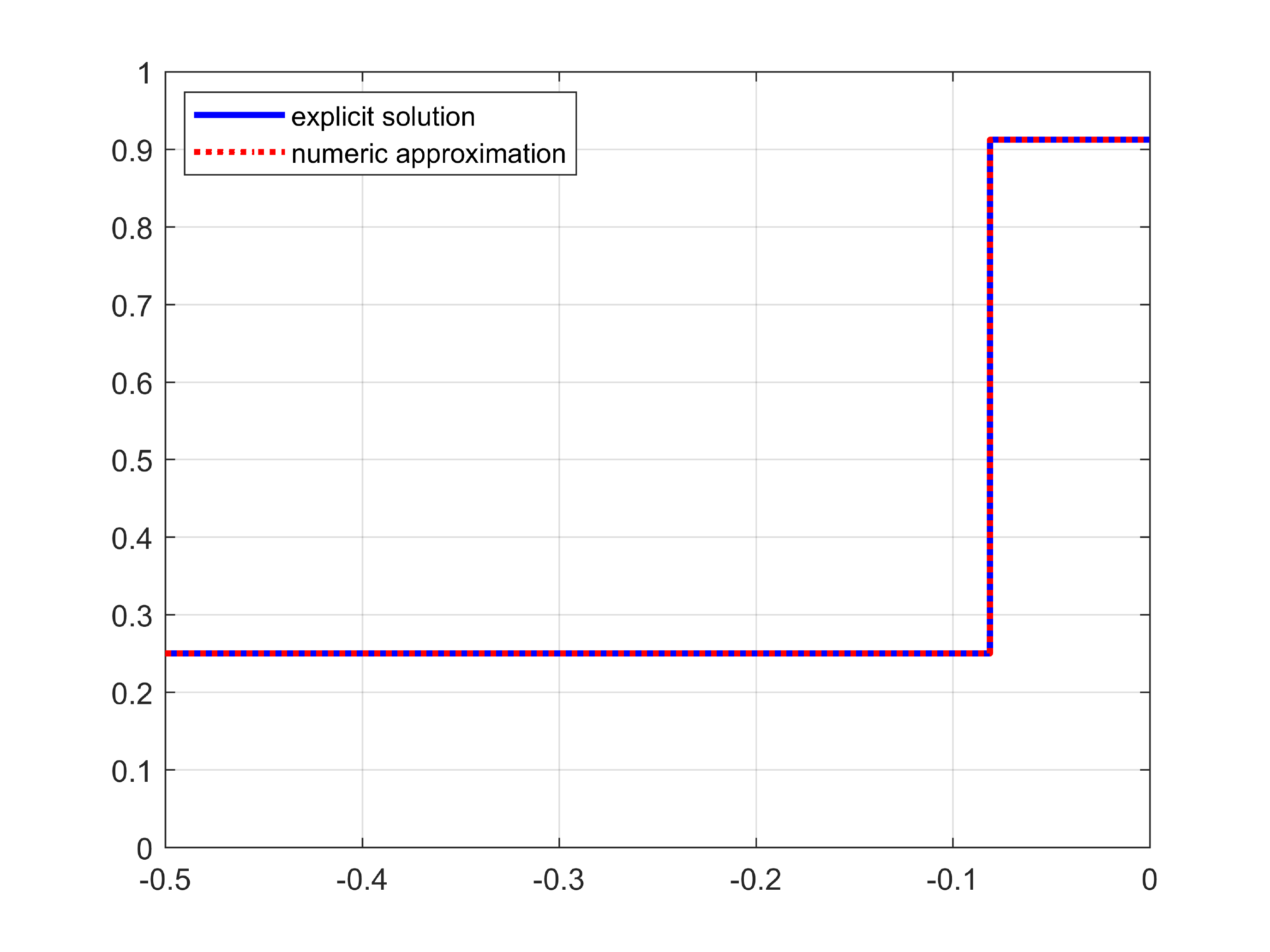}
	\end{subfigure}
	\begin{subfigure}[b]{0.328\textwidth}
			\includegraphics[width=\textwidth]{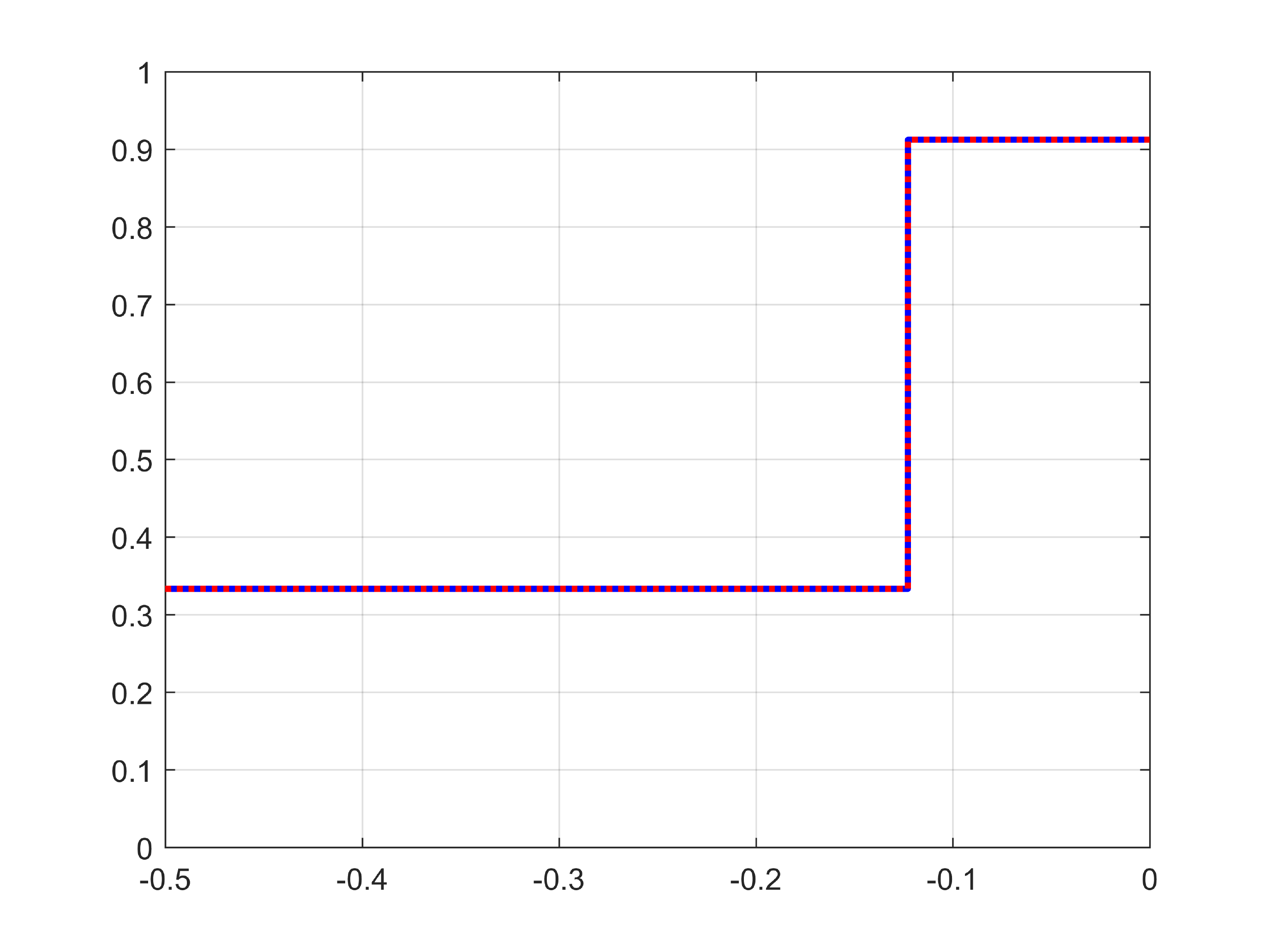}
	\end{subfigure}
	\begin{subfigure}[b]{0.328\textwidth}
			\includegraphics[width=\textwidth]{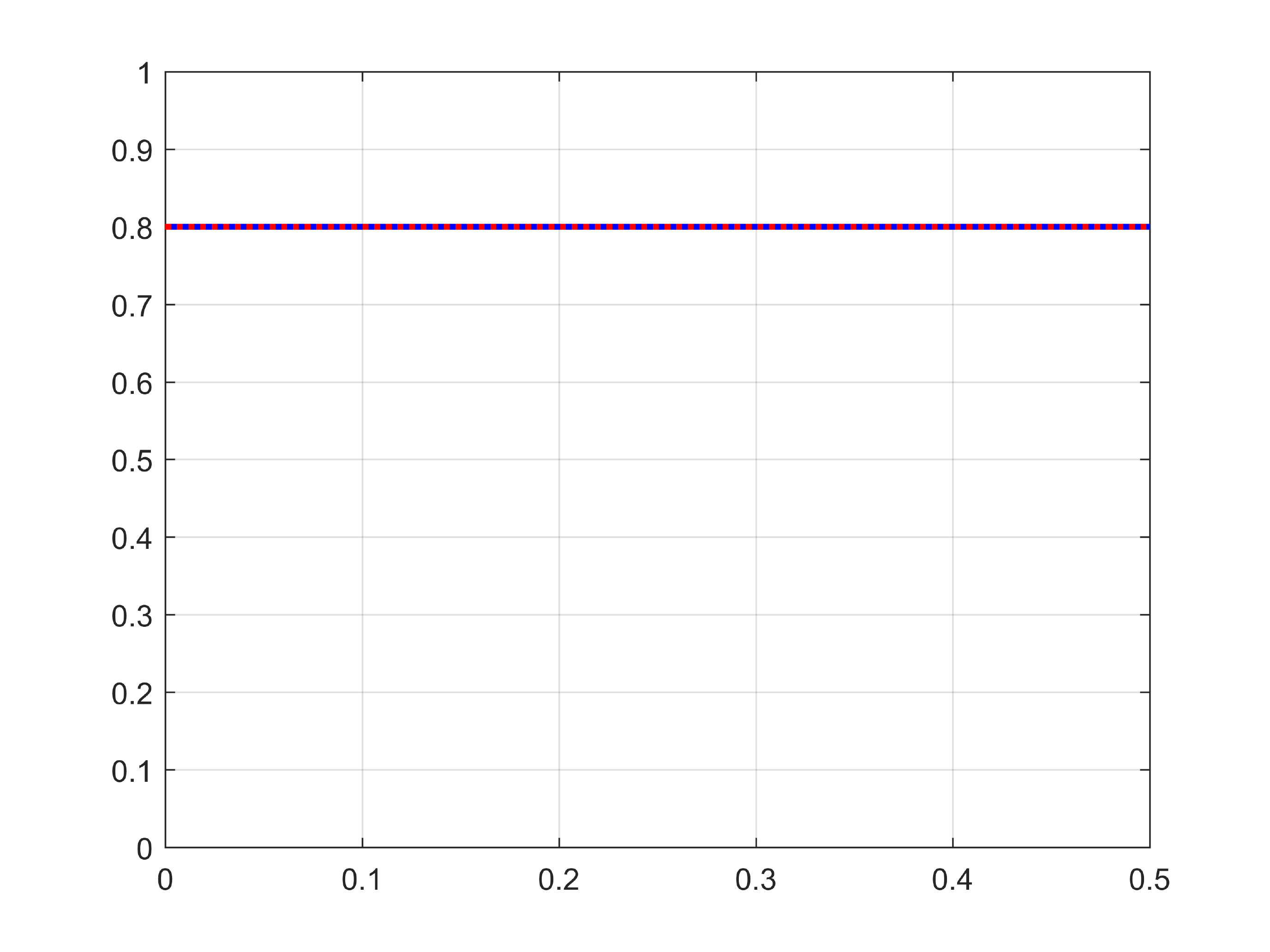}
	\end{subfigure}
	\\
	\begin{subfigure}[b]{0.328\textwidth}
			\includegraphics[width=\textwidth]{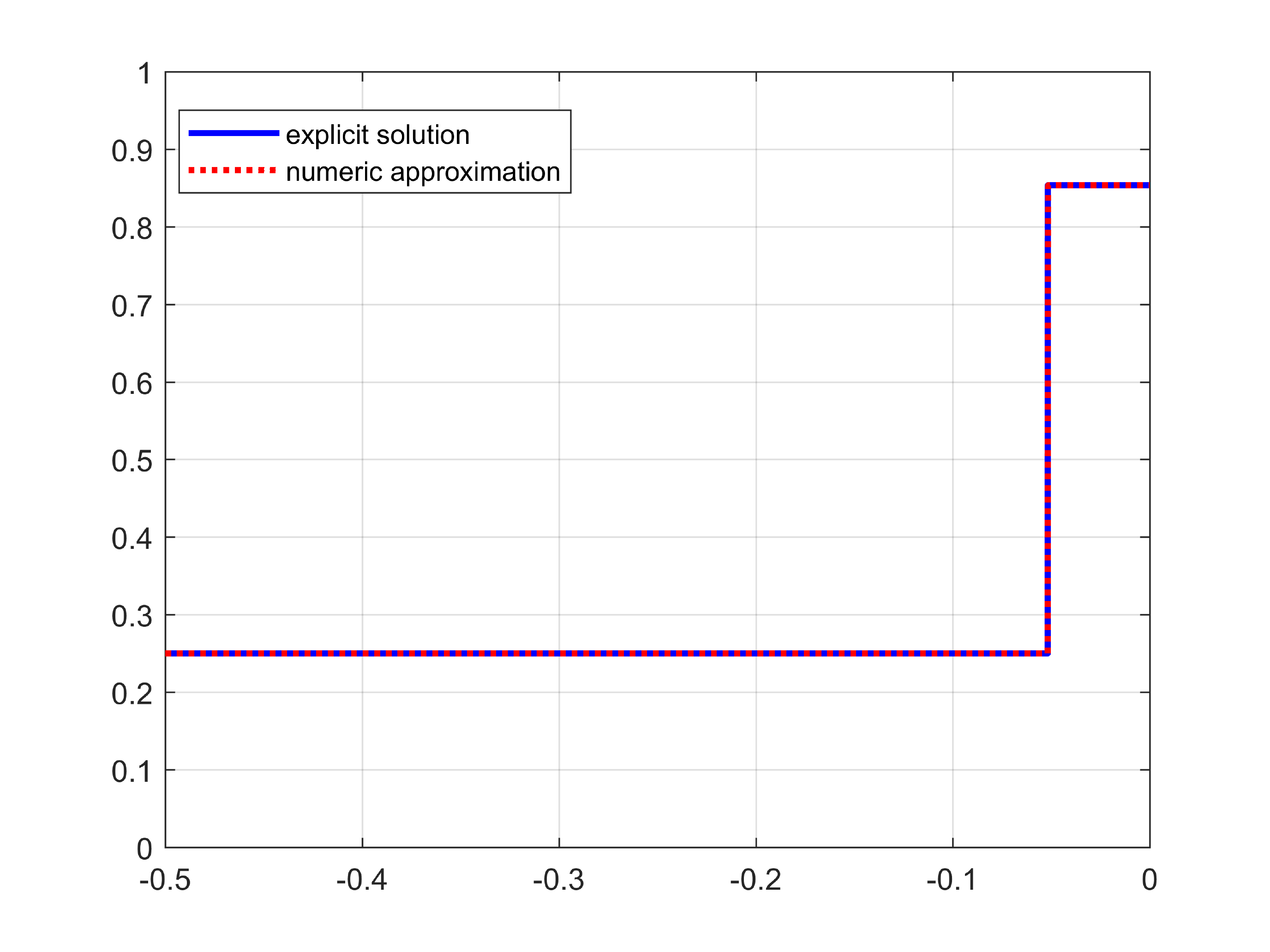}
	\end{subfigure}
	\begin{subfigure}[b]{0.328\textwidth}
			\includegraphics[width=\textwidth]{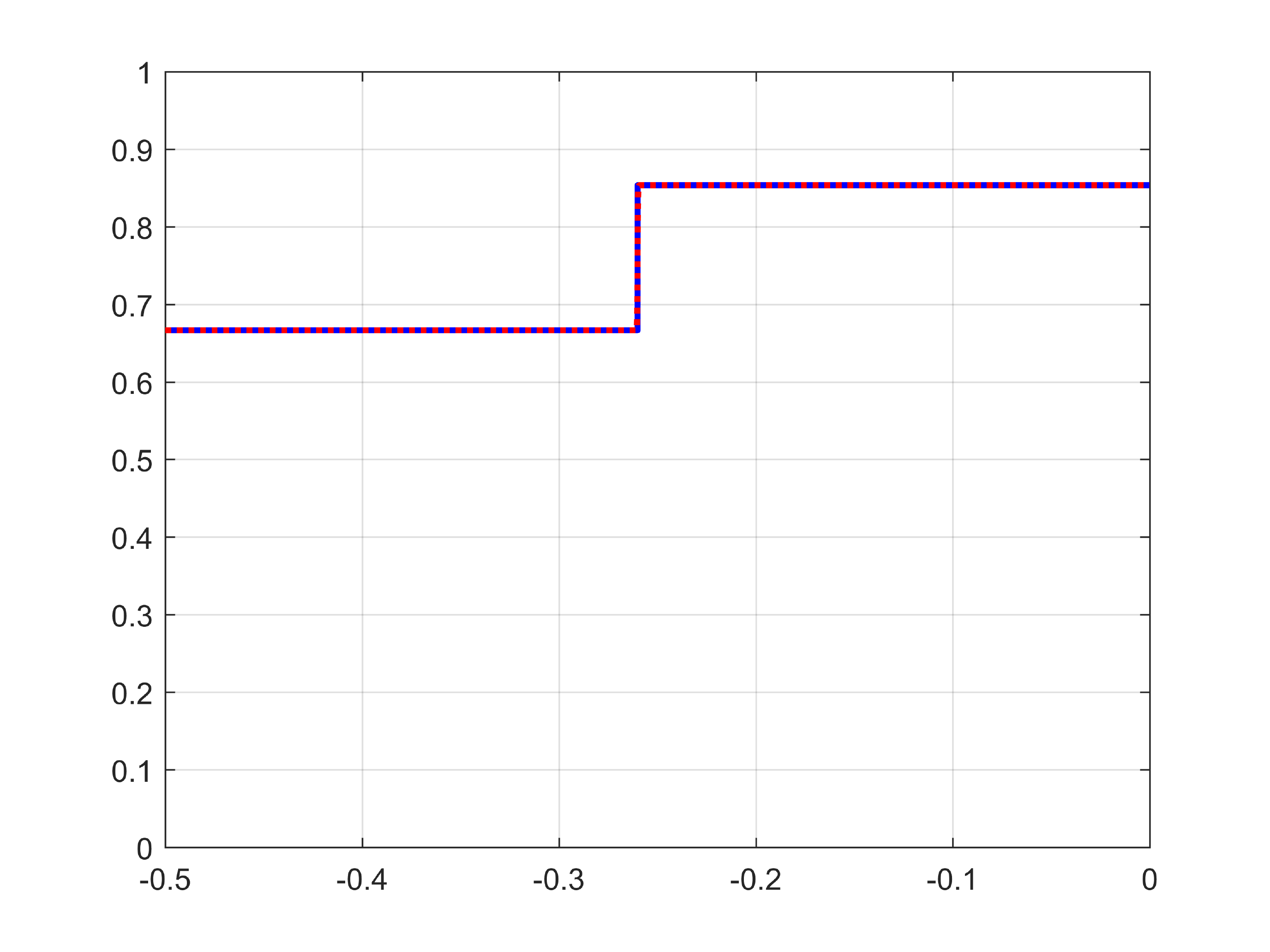}
	\end{subfigure}
	\begin{subfigure}[b]{0.328\textwidth}
			\includegraphics[width=\textwidth]{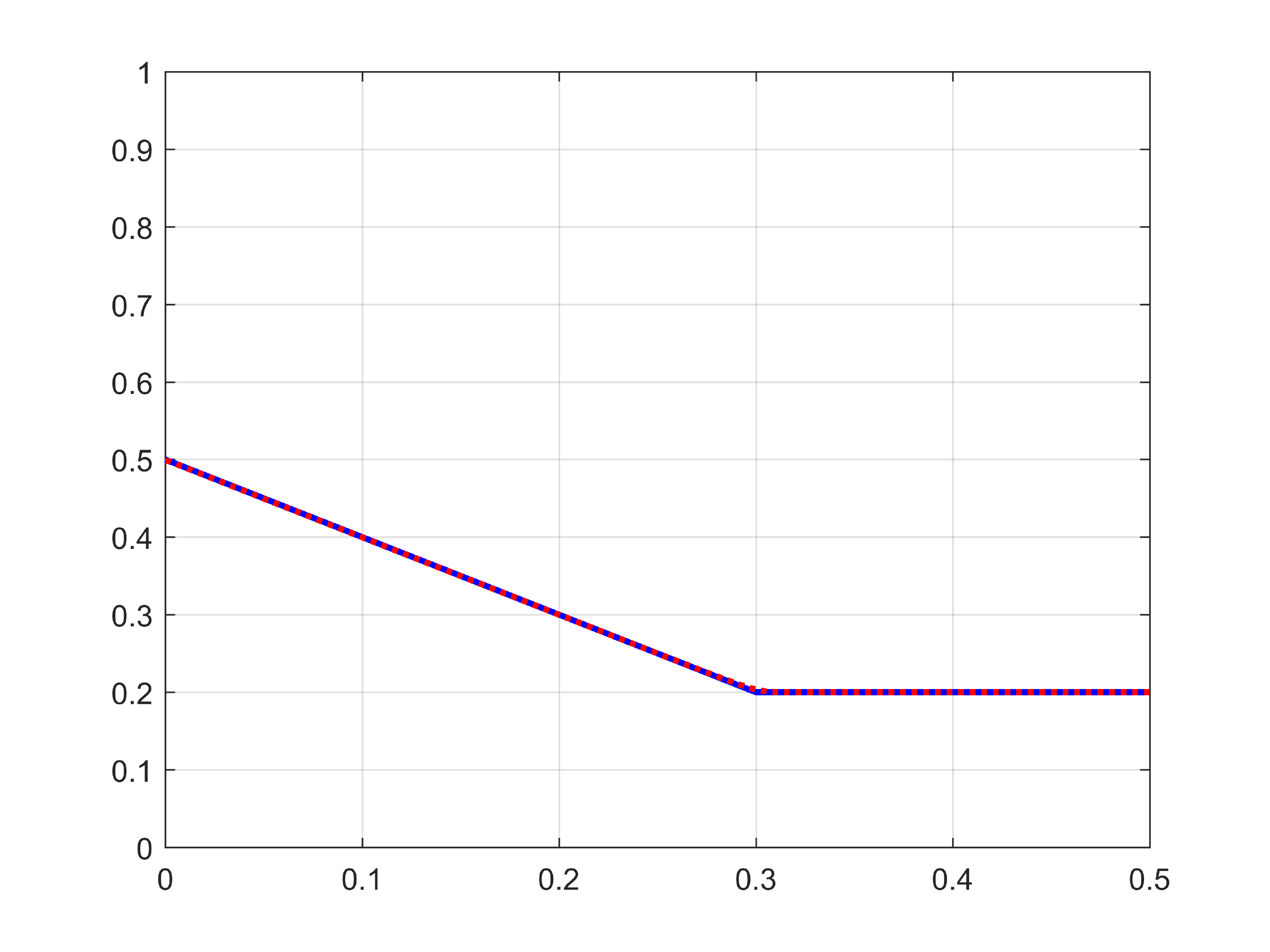}
	\end{subfigure}
	\caption[The comparison between the explicit solution and the numeric approximation for the Riemann problems at time $t=1/2$]{With reference to the Riemann problems test for a $2$-$1$ merge: comparison between the explicit solution and the numeric approximation for the Riemann problems at time $t=1/2$. First line shows the comparison of the profiles of solution along each arc corresponding to the initial condition $\vec\rho_{0,a}$; while last line refers to the initial datum $\vec\rho_{0,b}$. 
	}
	\label{fig:Riemann2-1}
\end{figure}

In Figures~\ref{fig:Riemann2-1} we present a qualitative comparison between the numerically computed solution and the explicitly one at time $t=\frac{1}{2}$, corresponding to the Riemann problems having $\vec\rho_{0,a} =\left(1/4,1/3,4/5\right)$ and $\vec\rho_{0,b}=\left(1/4,2/3,1/5\right)$ as initial conditions, respectively.
		%
		%
		%
		%
		%
	
We observe good agreement between the exact and the numeric solution. The parameters for the simulations are $\Delta x =10^{-4}$ and $\Delta t = 0.5\, \Delta x =0.5 \times 10^{-4}$.

\subsection{Riemann problem for a $1$-$2$ divide.}\label{sec:1-2}
We consider a network consisting of three edges and one junction, with one incoming and two outgoing arcs. We consider $[-\frac{1}{2}, 0]$ as domain of computation for the incoming arc and $[0,\frac{1}{2}]$ for the outgoing ones, and we take a normalized flux $f(\rho) = \rho \left(1 - \rho \right)$ for each arc. We consider as initial conditions for the Riemann problems $\vec\rho_{0,a}=\left(1/4,2/3,4/5\right)$ and $\vec\rho_{0,b} = \left(3/4,1/3,4/5\right)$, respectively.
	
		%
		%
		%
		%
		%
	%

\begin{figure}
	\centering
	\begin{subfigure}[b]{0.328\textwidth}
			\includegraphics[width=\textwidth]{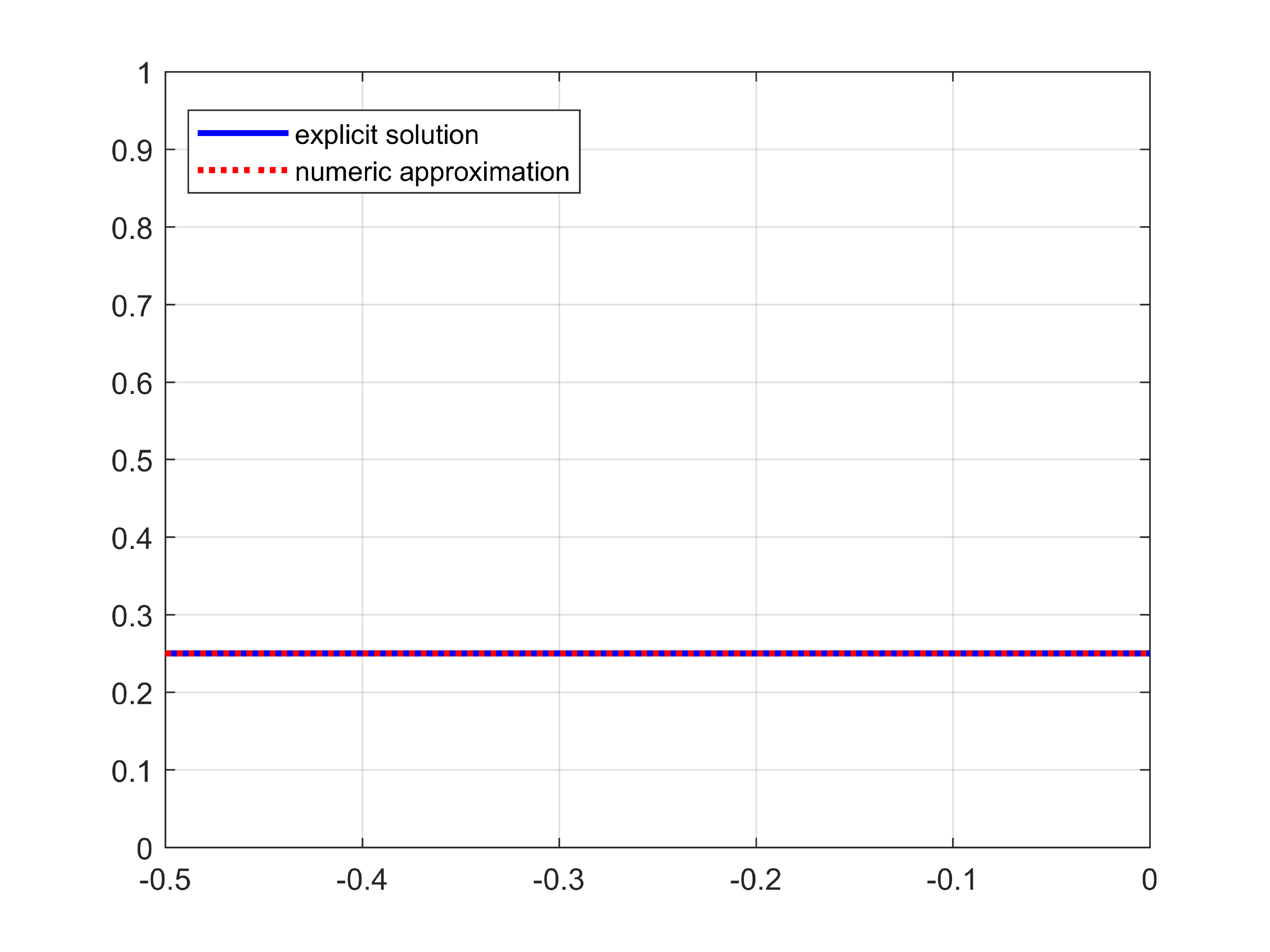}
	\end{subfigure}
	\begin{subfigure}[b]{0.328\textwidth}
			\includegraphics[width=\textwidth]{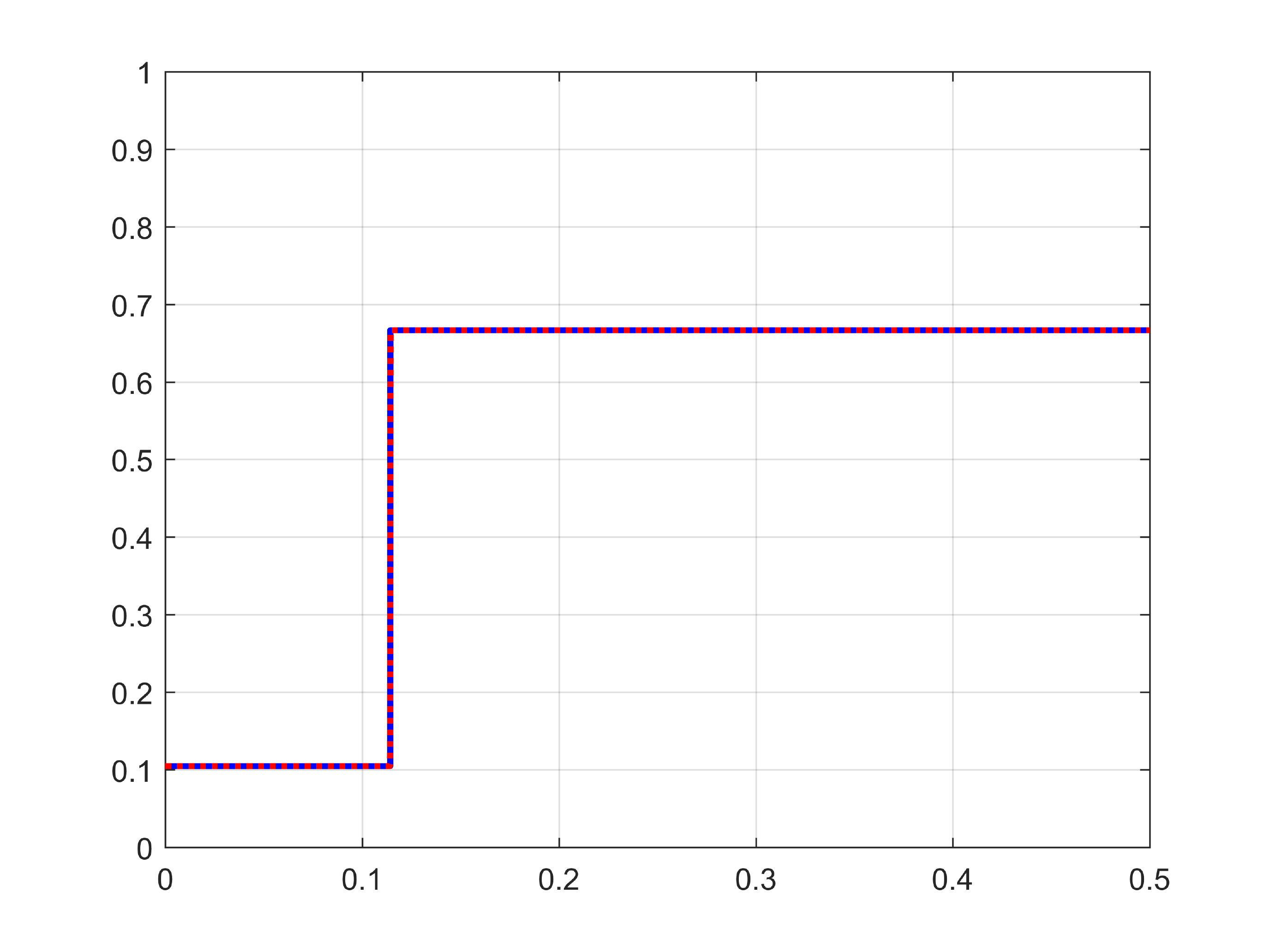}
	\end{subfigure}
	\begin{subfigure}[b]{0.328\textwidth}
			\includegraphics[width=\textwidth]{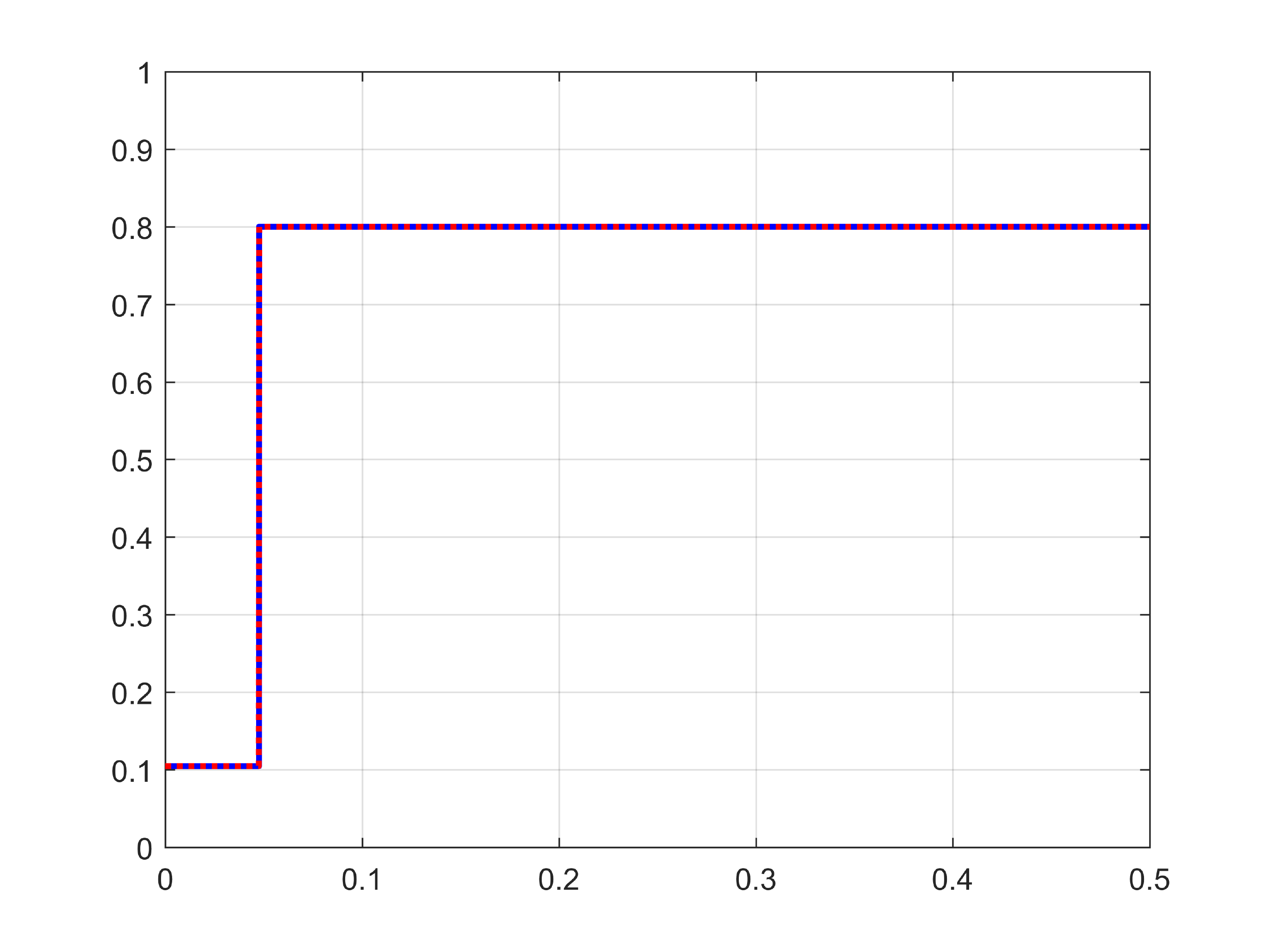}
	\end{subfigure}
	\\
	\begin{subfigure}[b]{0.328\textwidth}
			\includegraphics[width=\textwidth]{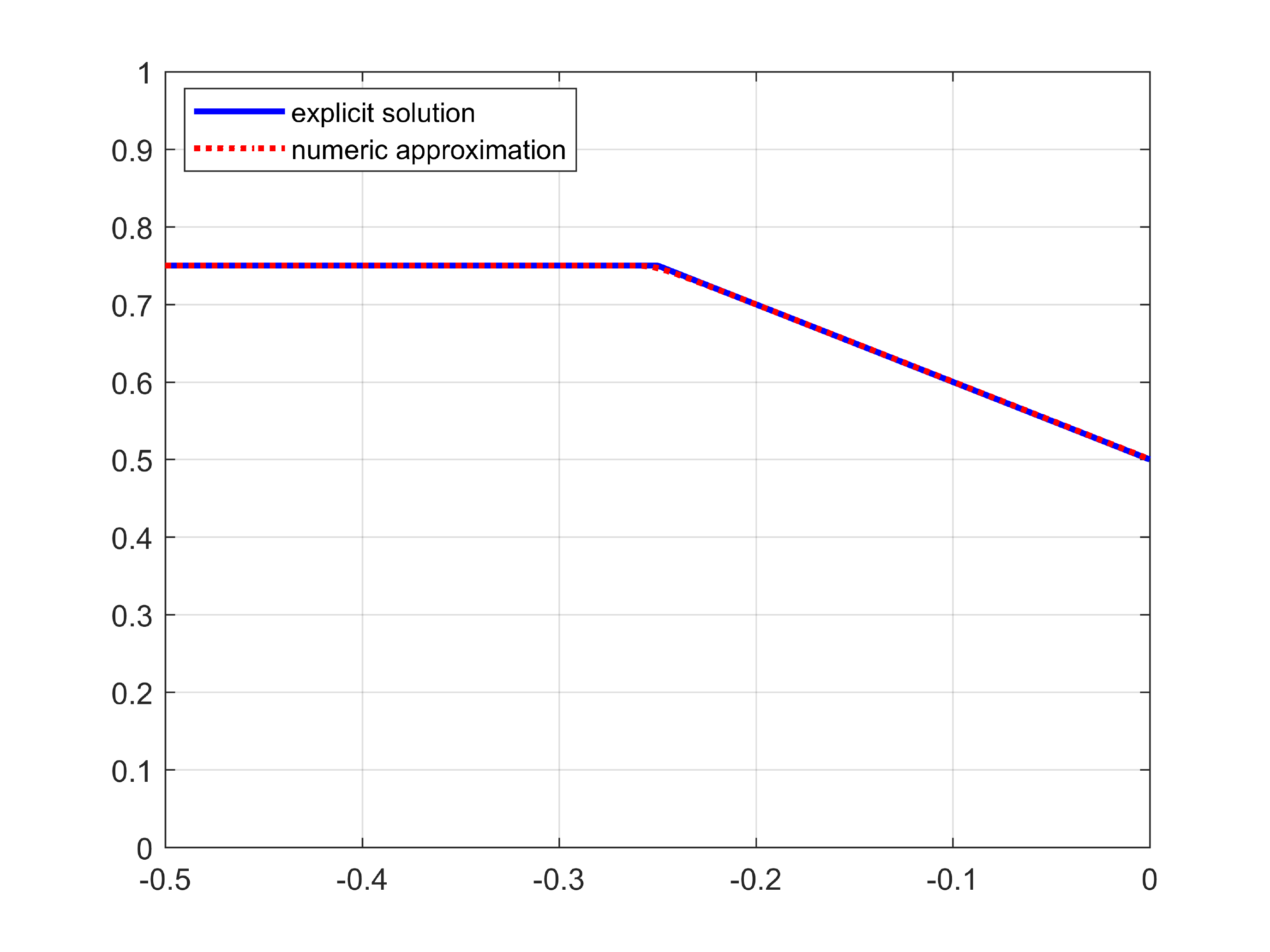}
	\end{subfigure}
	\begin{subfigure}[b]{0.328\textwidth}
			\includegraphics[width=\textwidth]{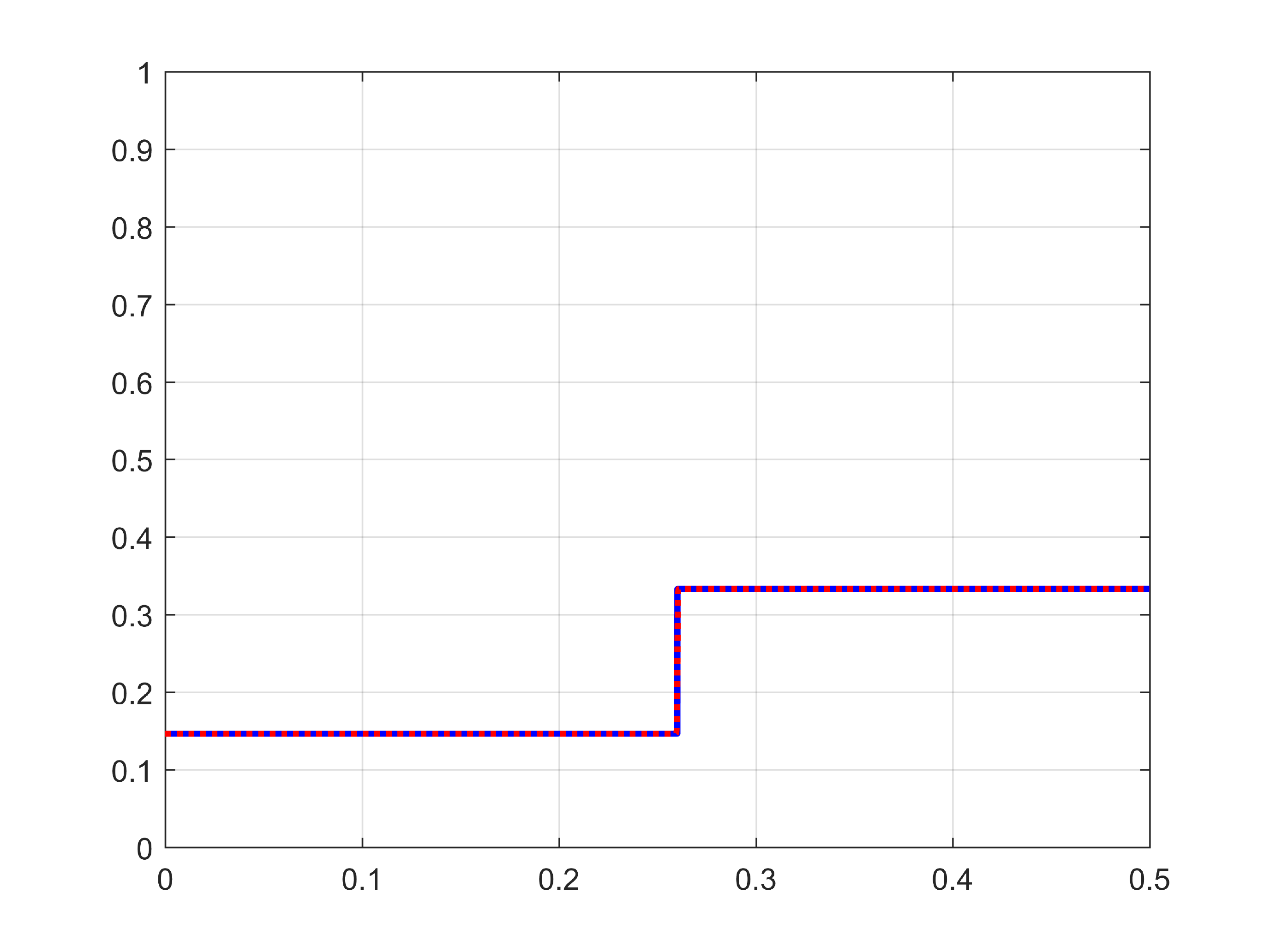}
	\end{subfigure}
	\begin{subfigure}[b]{0.328\textwidth}
			\includegraphics[width=\textwidth]{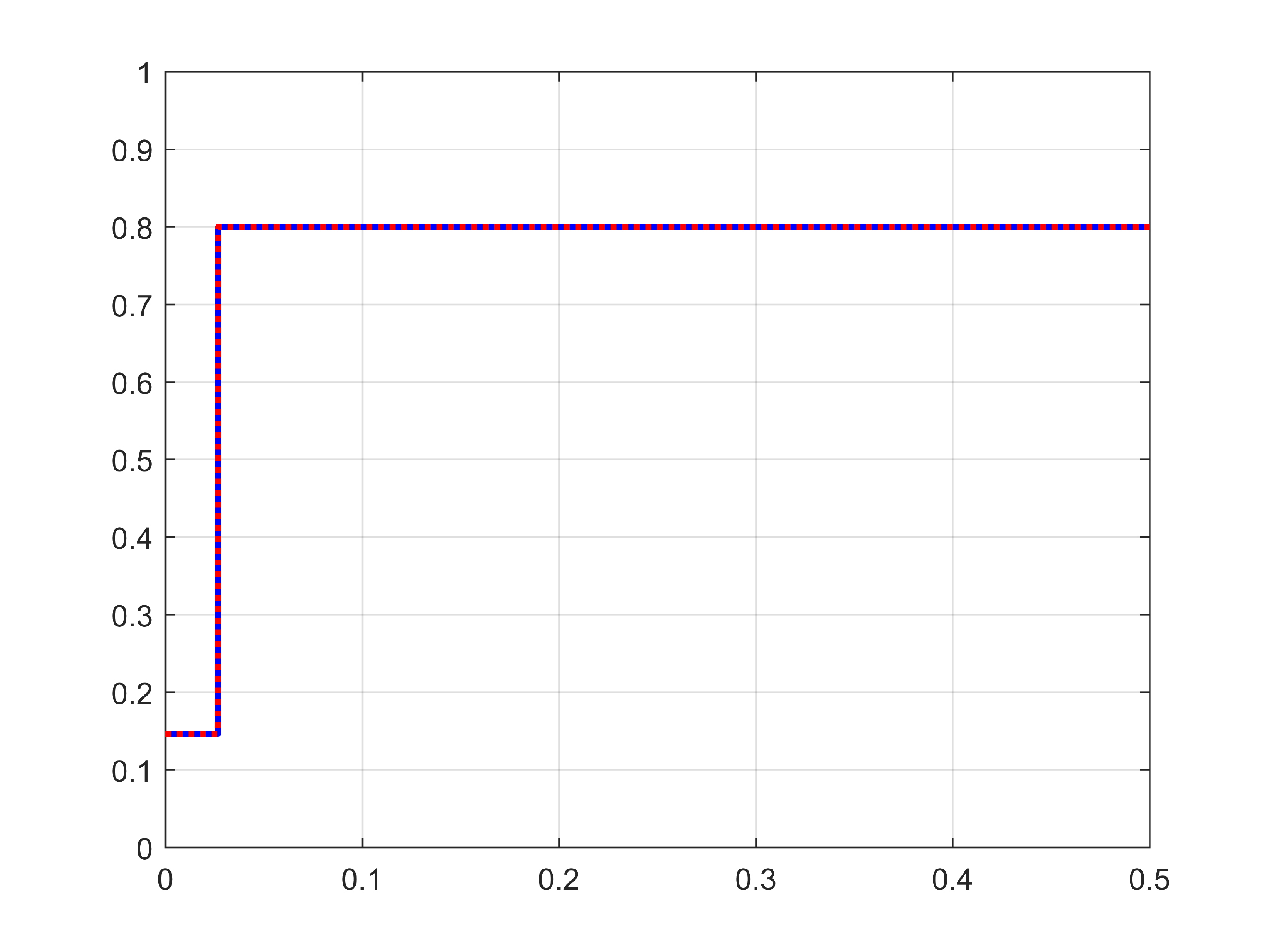}
	\end{subfigure}
	\caption[The comparison between the explicit solution and the numerical approximation for the Riemann problems at time $t=1/2$]{With reference to the Riemann problems for a $1$-$2$ divide: comparison between the explicit solution and the numerical approximation for the Riemann problems at time $t=1/2$. First line shows the comparison of the profiles of solution along
each arc corresponding to the initial condition $\vec\rho_{0,a}$; while last line refers to
the initial datum $\vec\rho_{0,b}$.
	}
	\label{fig:Riemann1-2}
\end{figure}
	
Figure \ref{fig:Riemann1-2} shows a qualitative comparison between the numerically computed solution and the explicitly one at time $t=\frac{1}{2}$, corresponding to the above initial conditions. Also in this case, we can observe good agreements between the exact solution and its numerical approximation. The parameters for the computed solution are $\Delta x =10^{-4}$ and $\Delta t = 0.5\, \Delta x =0.5 \times 10^{-4}$.

\subsection{Riemann problem for $2$-$2$ network.}

We consider here a junction with two incoming and two outgoing arcs. We consider $[-1/6, 0]$ as domain of computation for the incoming arcs and $[0,1/6]$ for the outgoing ones, and we take a normalized flux $f(\rho) = \rho \left(1 - \rho \right)$ for each arc. As initial condition for the Riemann problems on the network we choose $\vec\rho_{0,a} = \left(1/4,1/5,2/3,5/6\right)$ and $\vec\rho_{0,b} = \left(3/4,1/5,1/3,1/6\right)$.
	%
	%
	%
	%
	%
	%
	%
	%
%

\begin{figure}
	\centering
			\includegraphics[width=.61\textwidth]{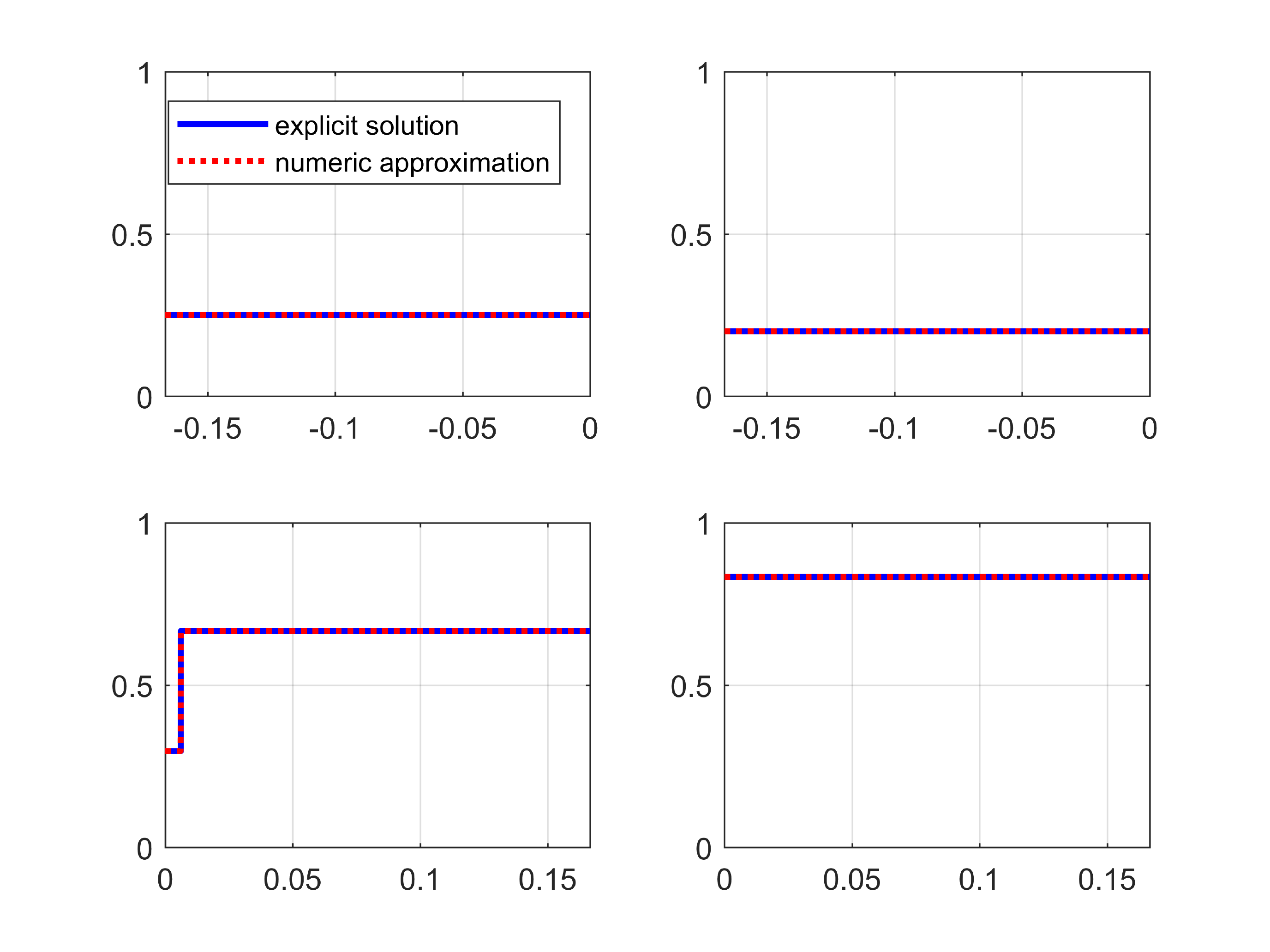}
	\caption[The comparison between the explicit solution and the numeric approximation for the Riemann problems at time $t=1/6$ corresponding to the initial condition $\vec\rho_{0,a}$]{With reference to the Riemann problems for a $2$-$2$ network: comparison between the explicit solution and the numeric approximation for the Riemann problems at time $t=1/6$ corresponding to the initial condition $\vec\rho_{0,a}$. First line shows the comparison of the profile of solution on the incoming arcs; last line displays the comparison of the profile of solution on the outgoing arcs.
	}
	\label{fig:Riemann2-2a}
\end{figure}

\begin{figure}
	\centering
			\includegraphics[width=.61\textwidth]{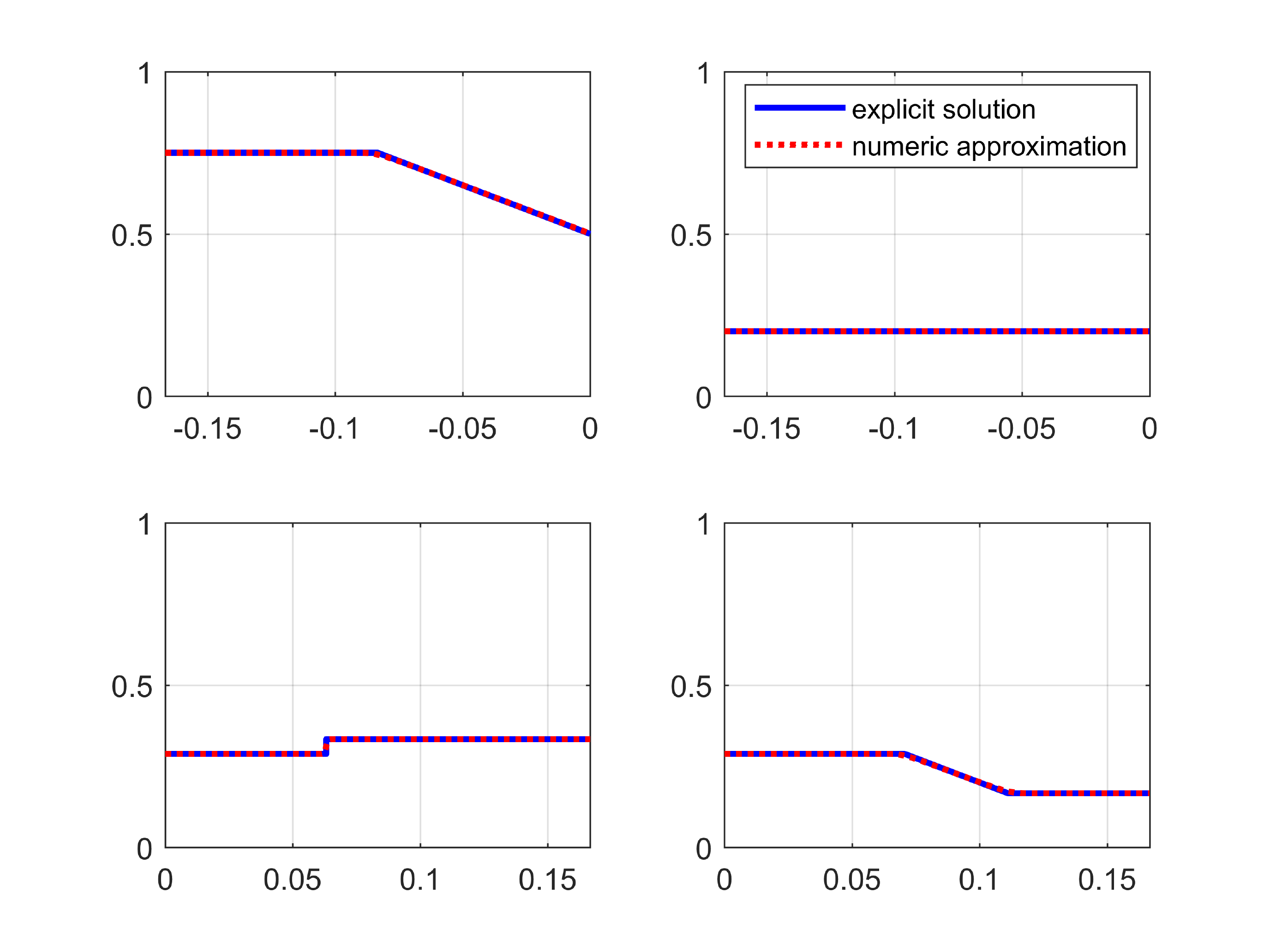}
	\caption[The comparison between the explicit solution and the numeric approximation for the Riemann problems at time $t=1/6$ corresponding to the initial condition $\vec\rho_{0,b}$]{With reference to the Riemann problems for a $2$-$2$ network: comparison between the explicit solution and the numeric approximation for the Riemann problems at time $t=1/6$ corresponding to the initial condition $\vec\rho_{0,b}$. First line shows the comparison of the profile of solution on the incoming arcs; last line displays the comparison of the profile of solution on the outgoing arcs.
	}
	\label{fig:Riemann2-2b}
\end{figure}

Also in this case, in Figures~\ref{fig:Riemann2-2a} and~\ref{fig:Riemann2-2b} we find a good agreement between the exact and the numerical solutions. The parameters for the simulation are $\Delta x =10^{-4}$ and $\Delta t = 0.5 \,\Delta x =0.5 \times 10^{-4}$.



\subsection{A network with no discontinuity at junction.}\label{sec:nojunction}

In this section we consider a network consisting of an arc without any discontinuity at the junction, namely, we assume that the flux on the incoming arc coincides with the flux on the outgoing arc, and therefore, this setting is equivalent to the case of a single arc. This simulation exploits the idea consisting in solving two scalar conservation laws on half-space coupled by an ad hoc transmission condition at the interface~\cite{AC}. We remark that, in the case of discontinuous flux, the transmission condition can be interpreted in terms of a flux constraint at the interface~\cite{AndreianovCances}.

We apply the scheme to two different domains, one including the junction located at $x=0$, and one not. More in details, we choose $[-1/2,1/2]$ and $[0,1]$ as domain of computation, and $f(\rho)=\rho(1-\rho)$ as flux along the arcs. We consider the initial densities $\rho_1^0=0.75\chi_{[-1/4,0]}$ and $\rho_2^0=0.75\chi_{[0,1/4]}$ for the first simulation, and $\rho_2^0=0.75\chi_{[0,1/4]}$ for the second one. The parameters of computation are $\Delta x=10^{-4}$ and $\Delta t=0.5\,\Delta x=0.5\times 10^{-4}$. 

\begin{figure}
\centering
	\begin{subfigure}[b]{0.327\textwidth}
			\includegraphics[width=\textwidth]{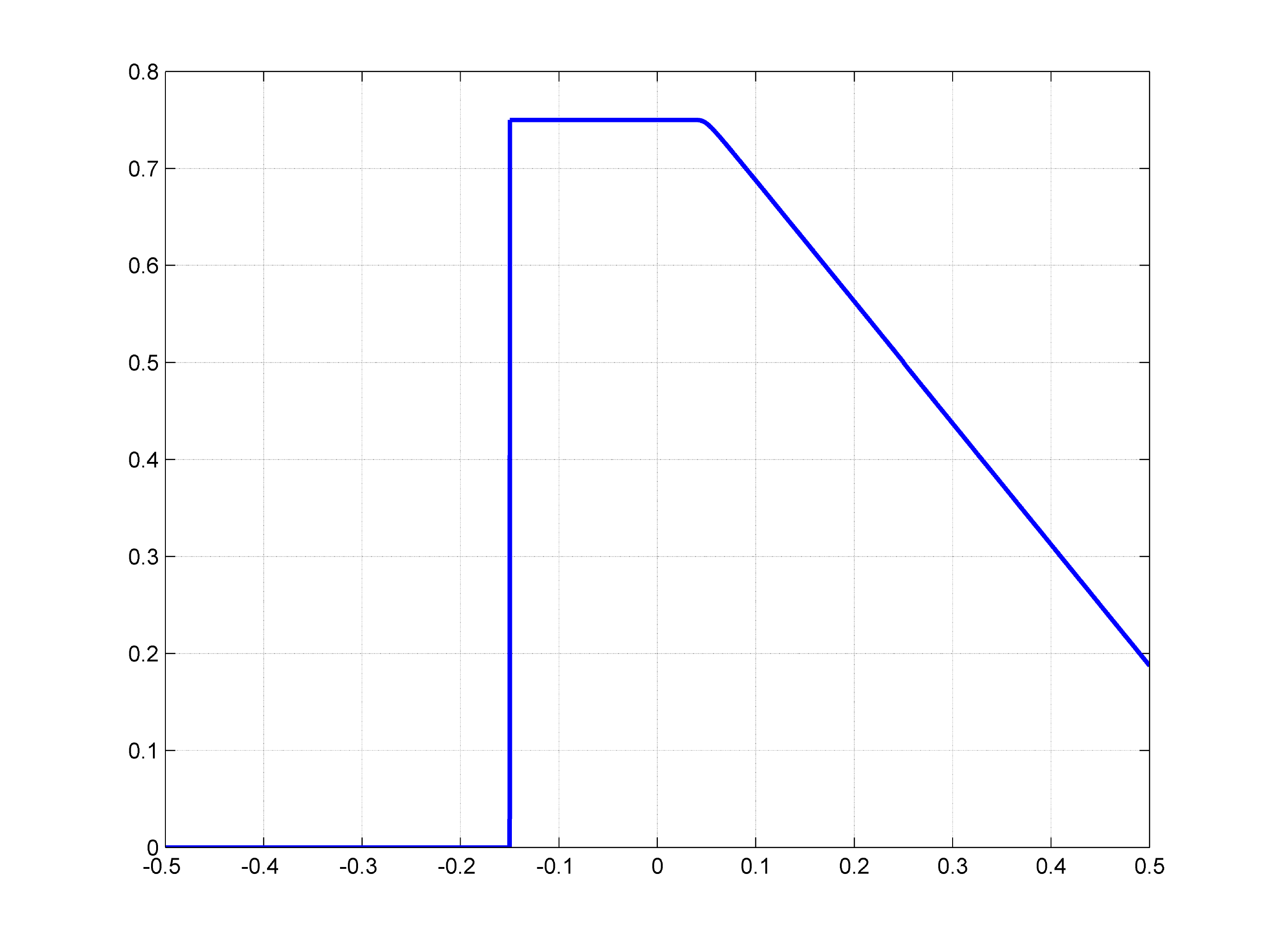}
			\caption{The profile of the solution of the first simulation at time $t=2/5$ for $x \in [-1/2,1/2]$.}
	\end{subfigure}
	\quad
	\begin{subfigure}[b]{0.327\textwidth}
			\includegraphics[width=\textwidth]{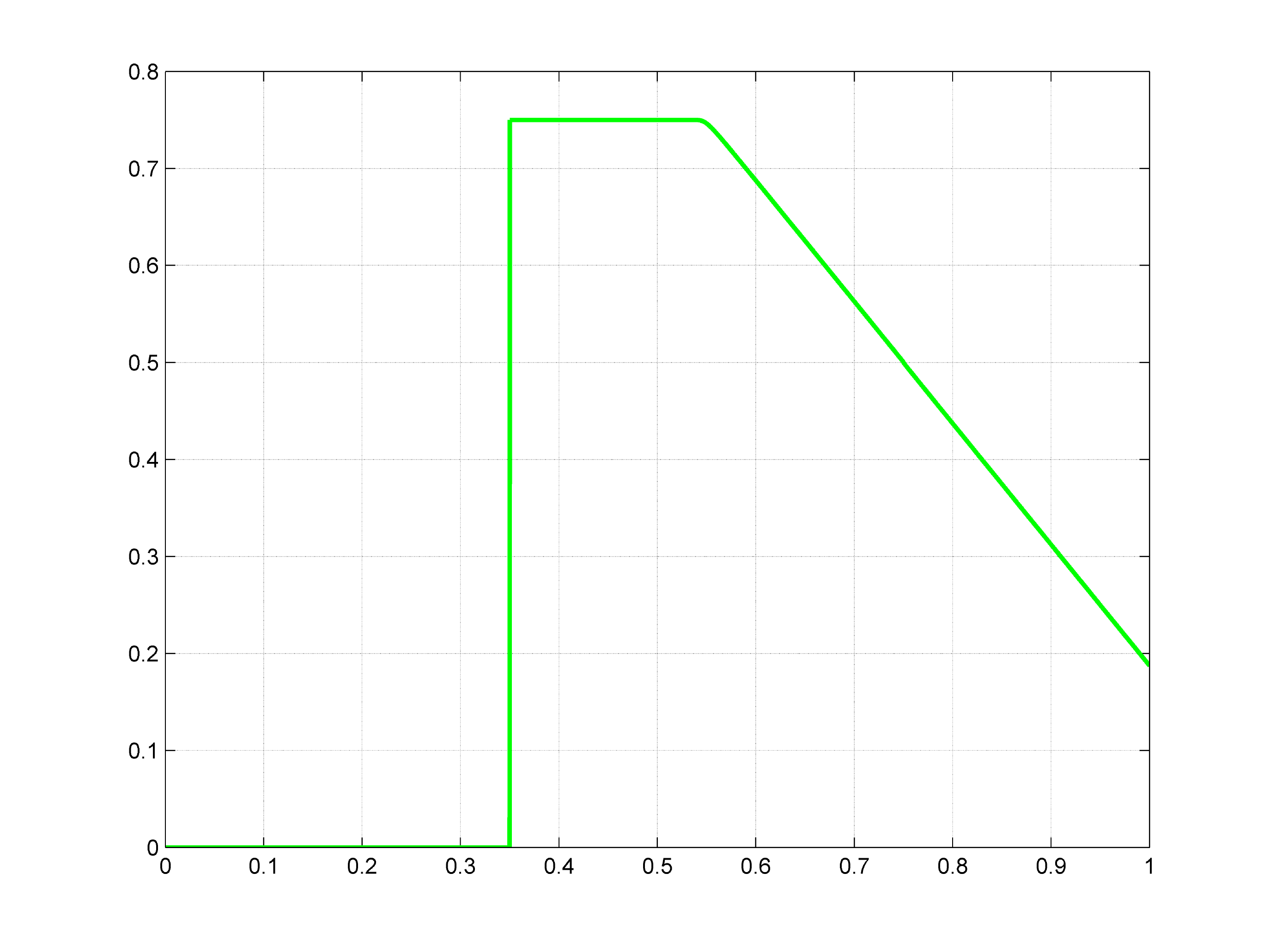}
			\caption{The profile of the solution of the second simulation at time $t=2/5$ for $x \in [0,1]$.}
	\end{subfigure}
	\\
	\begin{subfigure}[b]{0.327\textwidth}
			\includegraphics[width=\textwidth]{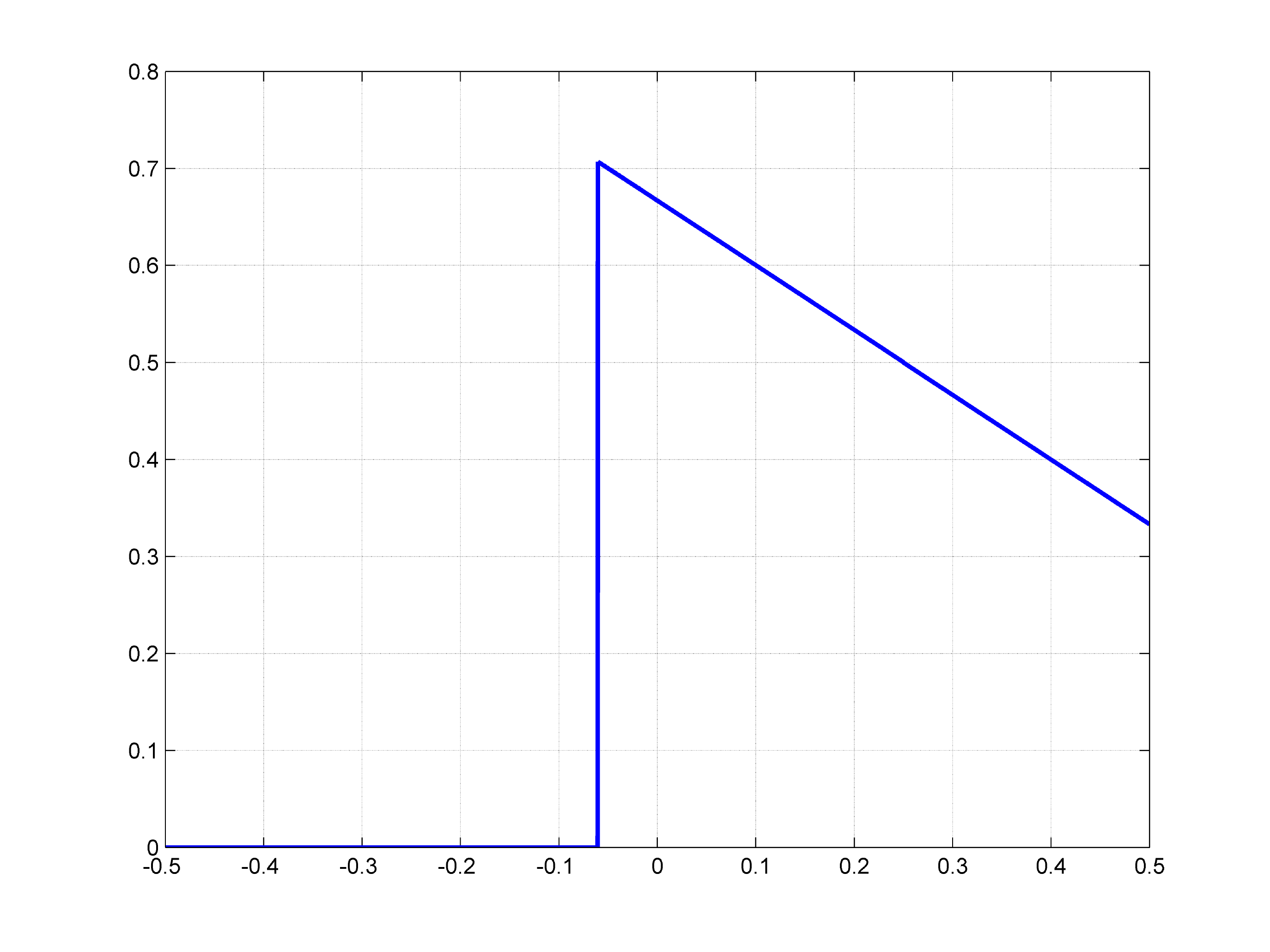}
			\caption{The profile of the solution of the first simulation at time $t=3/4$ for $x \in [-1/2,1/2]$.}
	\end{subfigure}
	\quad
	\begin{subfigure}[b]{0.327\textwidth}
			\includegraphics[width=\textwidth]{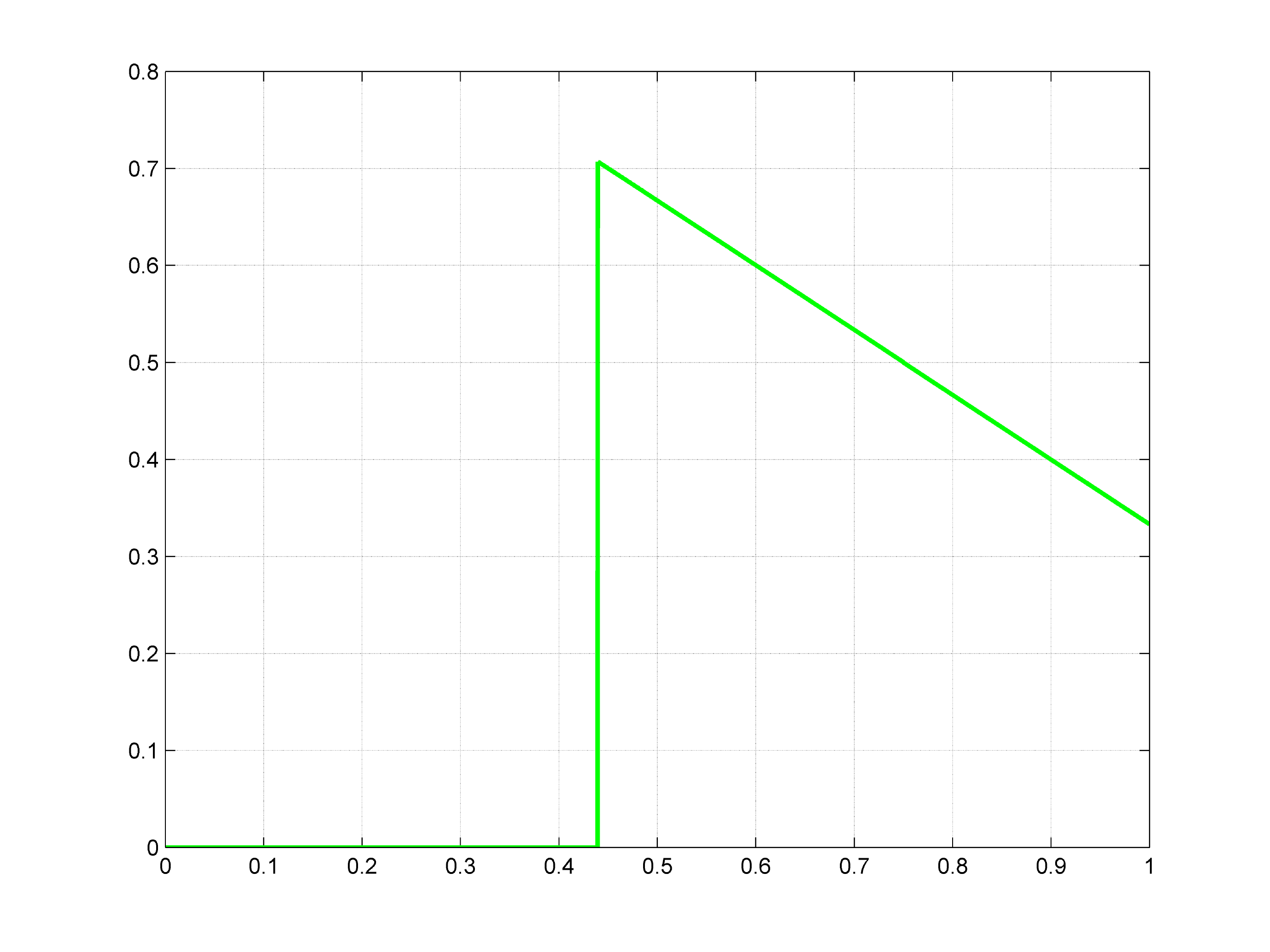}
			\caption{The profile of the solution of the second simulation at time $t=3/4$ for $x \in [0,1]$.}
	\end{subfigure}
	\caption[The case of a network with no discontinuity at the junction]{With reference to the case of a network with no discontinuity at the junction: comparison between the two simulations at two different times.}
	\label{fig:nojunc}
\end{figure}

Figure~\ref{fig:nojunc} displays the comparison between the profiles of solutions in the two simulations at times $t=2/5$ and $t=3/4$. We can observe the same qualitative behavior shifted of $|x|=1/2$.

\bibliographystyle{acm}
    \bibliography{validation_ACD}

\end{document}